\documentclass[12pt]{article}

\usepackage{notoccite}

\usepackage[longnamesfirst,sort]{natbib}
\bibpunct[, ]{(}{)}{;}{a}{,}{,}

\usepackage[title]{appendix}
\usepackage{tikz}
\usetikzlibrary{matrix}
\usepackage{calligra} 

\usepackage{multirow}
\usepackage{caption}
\usepackage{amsmath}
\usepackage{amsfonts}
\usepackage{amssymb}
\usepackage{graphicx}
\usepackage{graphics}
\usepackage{amscd,amsmath,amstext,amsfonts,amsbsy,amssymb,amsthm}
\usepackage[unicode,colorlinks,linkcolor= blue]{hyperref}
\usepackage{nccmath}
\usepackage{framed}
\usepackage{authblk}
\usepackage{bbm}
\usepackage{listings}
\usepackage{pdfpages}
\graphicspath{ {images/} }
\usepackage{tikz}
\usetikzlibrary{arrows,positioning,shapes.geometric}

\newtheorem{assumption}{Assumption}

\newtheorem{thm}{Theorem}[section]

\newtheorem{proposition}[thm]{Proposition}
\newtheorem{corollary}[thm]{Corollary}

\newcommand{\bequ}{\begin{equation}}
	\newcommand{\eequ}{\end{equation}}

\newcommand{\norm}[1]{\left\lVert#1\right\rVert}
\begin{document}

	
	\begin{center}
		{\Large
			{\sc  \bf{Aggregated kernel based tests for signal detection in a regression model}}
		}
		\bigskip
		
		Bui Thi Thien Trang $^{1}$ 
		\bigskip
		
		{\it
			$^{1}$ Institut de Math\'ematiques de Toulouse \\
			Universit\'e Paul Sabatier 118, route de Narbonne F-31062 Toulouse Cedex 9\\
			tbui@math.univ-toulouse.fr
		}
	\end{center}
	\bigskip
	
	
	{\bf Abstract.} Considering a regression model, we address the question of testing the nullity of the regression function. The testing procedure is available when the variance of the observations is unknown and does not depend on any prior information on the alternative. We first propose a single testing procedure based on a general symmetric kernel and an estimation of the variance of the observations. The corresponding critical values are constructed to obtain non asymptotic level-$\alpha$ tests. We then introduce an aggregation procedure to avoid the difficult choice of the kernel and of the parameters of the kernel. The multiple tests satisfy non-asymptotic properties and are adaptive in the minimax sense over several classes of regular alternatives.
	\smallskip
	
	

	\section{Introduction}\label{sec1}
	We observe $\left(X_i,Y_i\right)_{1\leq i \leq n}$ that obey to the regression model described as follows,
	\begin{equation}\label{Y1Y2-equation}
		Y_i= f\left(X_i\right)+\sigma \epsilon_i,\quad i=1,\cdots,n.
	\end{equation} 
	We assume that $X=\left(X_1,X_2,\cdots,X_n \right)$ are i.i.d real random variables with values in a measurable set $E$ such that $\left[0,1\right]\subset E \subset \mathbb{R}$ with bounded density $\nu$ with respect to the Lebesgue measure on $E$ and $\epsilon=\left(\epsilon_1,\epsilon_2,\cdots,\epsilon_n \right)$ are i.i.d standard Gaussian variables, independent of $\left(X_1,X_2,\cdots,X_n \right)$. All along the paper, $f$ is assumed to be in $\mathbb{L}^2\left(E,d\nu\right)$. We also assume that $\|f\|_\infty = \sup_{x\in E}|f(x)| < +\infty$. In order to estimate $\sigma^2$, we assume that we also observe $\left(Y^{'}_1,\cdots,Y^{'}_n \right)$ that obey to the model 
	\begin{equation}\label{y1y2-equation}
	Y_i^{'} = f\left(\frac{i}{n} \right) + \sigma \epsilon_i^{'} ,\quad i=1,\cdots,n,
	\end{equation} 
	where $\epsilon^{'}=\left(\epsilon_1^{'},\cdots,\epsilon_n^{'} \right)$ is independent of $\left(X_1,\cdots,X_n,\epsilon_1,\cdots,\epsilon_n \right)$. \\
	Given the observation of $\left(X_i,Y_i\right)_{1\leq i \leq n},\ \left(Y^{'}_i\right)_{1\leq i\leq n}$, we want to test the null hypothesis
	\begin{equation*}
		\left(H_0\right):\  f=0,
	\end{equation*}
	against the alternative
	\begin{equation*}
		\left(H_1\right):\  f\neq 0.
	\end{equation*}
	Hypothesis testing in nonparametric regression have been considered in the papers by \cite{king1988test}, \cite{hardle1990semiparametric}, \cite{hall1990bootstrap}, \cite{king1991testing} and \cite{delgado1992testing}. Tests for no effect in nonparametric regression are investigated in \cite{eubank1993testing}. In the paper of \cite{spokoiny1996adaptive}, the authors considered the particular case where $\sigma$ is assumed to be known. They propose tests that tests achieve the minimax rates of testing [up to an unvoidable $\log \log(n)$ factor] for a wide range of Besov classes. \cite{baraud2003adaptive} propose a test, based on model selection methods, for testing in a fixed design regression model that $\left(f\left(X_1\right),\cdots,f\left(X_n\right) \right)$ belongs to a linear subspace of $\mathbb{R}^n$ againts a nonparametric alternative. They obtain optimal rates of testing are up to a possible $\log n$ factor over various classes of alternatives simultaneously. More recently, in a Poisson process framework, \cite{fromont2012kernels, fromont2013two} consider two independent Poisson processes and address the question of testing equality of their respective intensities. They introduce tests based on a single kernel function and aggregate several kernel based tests to obtain adaptive minimax testing procedures over alternatives based on Besov or Sobolev balls.
	
	Our this work, we propose to construct aggregated kernel based testing procedures of $(H_0)$ versus $(H_1)$ in a regression model. Our test statistics are based on a single kernel function which can be chosen either as a projection or Gaussian kernel and we propose an estimation for the unknown variance $\sigma^2$. Our tests are exactly (and not only asymptotically) of level $\alpha$. We obtain the optimal non-asymptotic conditions on the alternative which guarantee that the probability of second kind error is at most equal to a precribed level $\beta$. However, the testing procedures that we introduce hereafter also intended to overcome the question of calibrating the choice of kernel and/or the parameters of the kernel. They are based on an aggregation approach, that is well-known in adaptive testing (\cite{baraud2003adaptive} and \cite{fromont2013two}). This paper is strengly inspired by the paper of \cite{fromont2013two}. Instead of considering a particular single kernel, we consider a collection of kernels and the corresponding collection of tests, each with an adapted level of significance. We then reject the null hypothesis when there exists at least one of the tests in the collection which rejects the null hypothesises. The aggregated testing procedures are constructed to be of level $\alpha$ and the loss in second kind error due to the aggregation, when unavoidable, is as small as possible. Then we prove that these multiples tests satisfy the adaptive minimax properties over several classes of alternatives. At last, we compare our tests with tests investigated in \cite{eubank1993testing} from a practical point of view.
	
	The paper is organized as follows. We describe the single tests based on a single kernel function with the corresponding critical values approximated by a Monte Carlo method in Section \ref{compactsect}. In Section \ref{sect3}, we specify the performances of the single tests for two particular examples of kernels and explain the reasons why we need to aggregate tests based on a collection of kernel functions which are presented in Section \ref{SectAggregateTest}. We present the simulation study in Section \ref{simulationSect} and the major proofs are given in Appendix. 
	\section{Single tests based on a single kernel.}\label{compactsect}
	
	\subsection{Definition of the testing procedure.}
	We assume that we observe $\{(X_i,Y_i)\}_{i=1}^n$ that obey to model (\ref{Y1Y2-equation}).	In order to estimate the unknown variance $\sigma^2$, we assume that we observe another sample $\left(Y^{'}_i \right)_{1\leq i\leq n}$ from the model (\ref{y1y2-equation}).
	We are interested in testing the null hypothesis $(H_0):\ f=0$ against $(H_1):\ f \neq 0$. Let $K$ be a symmetric kernel function: $E \times E\rightarrow \mathbb{R}$ satisfying:
	\begin{assumption}\label{assump1}
		\begin{equation*}
		\int_{E^2}^{} K^2(x,y)f(x)f(y)d\nu(x)d\nu(y) <+\infty.
		\end{equation*}
	\end{assumption}
	We introduce the test statistic $V_K$ defined as follows,
	\begin{equation}\label{Vk}
	V_K = \frac{T_K}{\hat{\sigma}^2_n},
	\end{equation} 
	where
	\begin{equation}\label{teststat}
		T_K = \frac{1}{n(n-1)} \sum_{i\neq j =1}^{n} K(X_i,X_j) Y_i Y_j
	\end{equation}
	and
	\begin{equation}\label{hat.sigma}
	\hat{\sigma}^2_n = \frac{1}{n}\sum_{i=1}^{n/2} \left(Y^{'}_{2i-1}-Y^{'}_{2i}\right)^2,
	\end{equation}
	where for the sake of simplicity, we assume that $n$ is even. Let us now introduce some notations. We set $K_{ij}=K\left(X_i,X_j\right),\ f_i=f(X_i)$ and $C(a,b)$ is a constant depending on $a$ and $b$, that will be used all along the paper and may vary from line to line. \\
	The expectation of $T_K$ is equal to
	\begin{align*}\label{varep}
		\mathbb{E}\left[T_K\right] &=\mathbb{E} \left[\mathbb{E}\left[\frac{1}{n(n-1)}\sum_{i \neq j=1}^{n}K_{ij} \left(f(X_i)+\sigma \epsilon_i \right) \left(f(X_j)+\sigma \epsilon_j \right) \bigg | X\right]\right]\\
		{} &=\mathbb{E}\left[\frac{1}{n(n-1)}\sum_{i \neq j=1}^{n}K_{ij} f(X_i)f(X_j) \right]\\
		{} &= \int_{E^2}^{}K(x,y)f(x)f(y)d\nu(x)d\nu(y).
	\end{align*}
	In the following, we denote for all $x\in E$,
	\begin{equation*}
	K[f](x) = \int_{E}^{}K(x,y)f(y)d\nu(y),
	\end{equation*}
	and for all $f,g \in \mathbb{L}^2(E,d\nu)$ 
	\begin{equation*}
	\langle f,g \rangle = \int_{E}^{}f(x)g(x)d\nu(x) \ \text{and}\  \|f\|^2=\langle f,f\rangle.
	\end{equation*}
	Within these notations, 
	\begin{equation}\label{var.ep}
	\mathbb{E}\left(T_K\right) = \langle K[f],f \rangle,
	\end{equation}
	whose existence is ensured by Assumption \ref{assump1}. We now compute the expectation of $\widehat{\sigma}^2_n$.
	\begin{align*}
	\mathbb{E}\left[\widehat{\sigma}^2_n\right] &= \mathbb{E}\left[\frac{1}{n}\sum_{i=1}^{n/2}\left[\left( f\left(\frac{2i-1}{n} \right) +\sigma \epsilon^{'}_{2i-1} \right) -\left( f\left(\frac{2i}{n} \right)+\sigma \epsilon^{'}_{2i} \right) \right]^2 \right]\\
	{} &= \frac{1}{n} \sum_{i=1}^{n/2}\left[ f\left(\frac{2i-1}{n} \right) -  f\left(\frac{2i}{n} \right) \right]^2 + \frac{1}{n}\sum_{i=1}^{n/2}\sigma^2\mathbb{E}(\epsilon^{'}_{2i-1}-\epsilon^{'}_{2i})^2\\
	{} &= a^2 +\sigma^2,
	\end{align*}
	with $a^2 := \frac{1}{n}\sum_{i=1}^{n/2}\left[f\left(\frac{2i-1}{n} \right) -f\left(\frac{2i}{n} \right) \right]^2$.\\
	Thus $\widehat{\sigma}^2_n$ is a biased estimator of $\sigma^2$ with bias $a^2$. If $f$ is a regular function this bias will be small.\\
	We have chosen to consider and study in this paper two possible examples of kernel functions. For each example, we give a simpler expression of $\mathbb{E}\left(T_K\right)$. \\
	
	\textit{Example 1}. When $E=[0,1]$, our first choice for $K$ is a symmetric kernel function based on a finite orthonormal family $\{\phi_\lambda,\ \lambda \in \Lambda \}$ with respect to the scalar product $\langle .,. \rangle$,
	\begin{equation}\label{projectionkernel}
		K(x,y)=\sum_{\lambda \in \Lambda}^{} \phi_\lambda (x)\phi_\lambda(y).
	\end{equation}
	For all $f$ in $\mathbb{L}^2([0,1],d\nu)$ we get
	\begin{align*}
	K[f](x) &= \int_{0}^{1} \left(\sum_{\lambda \in \Lambda}^{}\phi_\lambda(x)\phi_\lambda(y) \right)f(y)d\nu(y)\\
	{} &= \sum_{\lambda \in \Lambda}^{} \left(\int_{0}^{1}\phi_\lambda(y)f(y)d\nu(y) \right)\phi_\lambda(x) = \Pi_S(f),
	\end{align*}
	where $S$ is the subspace of $\mathbb{L}^2([0,1],d\nu)$ generated by the functions $\{\phi_\lambda,\ \lambda \in \Lambda \}$ and $\Pi_S$ denotes the orthogonal projection onto $S$ for $\langle . , . \rangle$. Thus
	\begin{equation*}
	\mathbb{E}\left(T_K\right) = \langle \Pi_S(f), f\rangle.
	\end{equation*}
	Hence, when $\{\phi_\lambda,\ \lambda \in \Lambda \}$ is well-chosen, $T_K$ can also be viewed as a relevant estimator of $\norm{f}^2$.\\
	
	\textit{Example 2.} When $E=\mathbb{R}$ and $\nu$ is a density function respect to the Lebesgue measure on $\mathbb{R}$, our second choice for $K$ is a Gaussian kernel defined by,
	\begin{equation}\label{13}
	K(x,y)=\frac{1}{h}k \left(\frac{x-y}{h}\right),\ \text{for all} \left(x,y\right)\in \mathbb{R}^2
	\end{equation}
	where $k(u)=\frac{1}{\sqrt{2\pi}}\exp \left(-u^2/2 \right), \ \text{for all}\ u\in \mathbb{R}$ and $h$ is a positive bandwidth. Then, for all $f\in \mathbb{L}^2(\mathbb{R},d\nu)$ we have
	\begin{equation*}
	K[f](x) = \int_{-\infty}^{\infty} \frac{1}{h} k\left(\frac{x-y}{h} \right)f(y)d\nu(y) = k_h\ast f(x),
	\end{equation*}
	where $\ast$ is the convolution producer with respect to the measure $\nu$ and $k_h(u)=\frac{1}{h}k\left(\frac{u}{h}\right),\ \forall u \in \mathbb{R}$. Thus in this case
	\begin{equation*}
	\mathbb{E}\left(T_K\right) = \langle k_h\ast f,f \rangle.
	\end{equation*}
	Hence, when the bandwidth $h$ is well chosen, $T_K$ can also be viewed as a relevant estimator of $\|f\|^2$.\\
	
	From the choices of the two examples above for $K$, we have seen that the test statistic $V_K$ can be viewed as a relevant estimator of $\|f\|^2$. Thus, it seems to be reasonable proposal to consider a test which rejects $(H_0)$ when $V_K$ is as "large enough". Now, we define the critical values used in our tests.
	
	We define
	\begin{equation}\label{test.epsilon}
	V_K^{\left(0\right)} =  \frac{\frac{1}{n(n-1)} \sum_{i\neq j =1}^{n} K(X_i,X_j) \epsilon_i \epsilon_j}{\frac{1}{n}\sum_{i=1}^{n/2} \left(\epsilon^{'}_{2i-1}-\epsilon^{'}_{2i}\right)^2}.
	\end{equation}
	Note that, under $\left(H_0 \right)$, conditionally on $X$, $V_K$ and $V_K^{\left(0\right)}$  have exactly the same distribution. We now choose the quantile of the conditional distribution of $V_K^{\left(0\right)}$ given $X$ as the critical value for our test. This quantity can easily be estimated by simulations.
	
	More precisely, for $\alpha$ in $(0,1)$, if $q^{(X)}_{K,1-\alpha}$ denotes the $(1-\alpha)$ quantile of the distribution of  $V_K^{\left(0\right)}$ conditionally on $X$, we consider the test that rejects $(H_0)$ when $V_K > q^{(X)}_{K,1-\alpha}$. The corresponding test function is defined by 
	\begin{equation}\label{testfunc}
		\Phi_{{K,\alpha}} =  \mathbbm{1} \{V_K > q^{(X)}_{K,1-\alpha}\}.
	\end{equation}
	Notice that in practice, the true quantile $q^{(X)}_{K,1-\alpha}$ is not available, but it can be approximated by a Monte Carlo procedure.
	
	\subsection{Probabilities of first and second kind errors of the test.}
	Since under $(H_0)$, $V_K$ and $V_K^{\left(0\right)}$ have the same distribution conditionally on $X$, for any $\alpha\in (0,1)$, we have
	\begin{equation*}
	\mathbb{P}_{(H_0)}\left(V_K > q^{(X)}_{K,1-\alpha} \bigg | X \right)\leq \alpha.
	\end{equation*}
	By taking the expectation over $X$, we obtain
	\begin{equation*}
	\mathbb{P}_{(H_0)}\left(\Phi_{K,\alpha} = 1\right) \leq \alpha.
	\end{equation*}

	Let us now consider an alternative hypothesis, corresponding to a non zero regression function $f$. Given $\beta$ in $(0,1)$, we now aim to determine a non-asymptotic condition on the regression function $f$ which guarantees that $\mathbb{P}_{f}(\Phi_{K,\alpha} = 0) \leq \beta$. Denoting by $q^\alpha_{K,1-\beta/2}$ the $(1-\beta/2)$ quantile of the conditional quantile $q^{(X)}_{K,1-\alpha}$,
	\begin{align*}
	\mathbb{P}_{f}\left(\Phi_{K,\alpha}=0\right) &= \mathbb{P}_{f}\left(V_K  \leq q^\alpha_{K,1-\beta/2} \right) + \mathbb{P}_{f}\left(V_K \leq q^{(X)}_{K,1-\alpha}, q^{(X)}_{K,1-\alpha} > q^\alpha_{K,1-\beta/2} \right)\\
	{} &\leq \mathbb{P}_{f}\left(V_K \leq q^\alpha_{K,1-\beta/2} \right) +\beta/2.
	\end{align*}
	Thus, a condition which guarantees that $\mathbb{P}_{f}\left(V_K \leq q^\alpha_{K,1-\beta/2} \right)\leq \beta/2$ will ensure that $\mathbb{P}_{f}(\Phi_{K,\alpha} = 0) \leq \beta$. The following proposition gives such a condition.
	
	\begin{proposition}\label{2nderror}
		Let $\alpha,\ \beta$ be the fixed levels in $(0,1)$. We have that
		\begin{equation*}
		\mathbb{P}_{f}\left(V_K \leq q^{\alpha}_{K,1-\beta/2} \right) \leq \beta/2,
		\end{equation*}
		as soon as
		\begin{equation}\label{12}
		\langle K\left[f\right],f \rangle \geq \sqrt{\frac{16A_K + 8B_K}{\beta}}+ D_{n,\beta}\  q^\alpha_{K,1-\beta/2} ,
		\end{equation}
		with
		\begin{align*}
		A_K &= \frac{n-2}{n(n-1)}\int_{E}^{}\left(K[f](x) \right)^2\left[f^2(x)+\sigma^2 \right]d\nu(x),\\
		B_K &= \frac{1}{n(n-1)}\int_{E^2}^{}K^2(x,y)\left[f^2(x)+\sigma^2 \right]\left[f^2(y)+\sigma^2 \right]d\nu(x) d\nu(y),\\
		D_{n,\beta} &= \sigma^2+a^2+\frac{4\sigma^2}{n}\sqrt{\left(\frac{n}{2}+\frac{na^2}{\sigma^2} \right)\ln \left(\frac{2}{\beta}\right)} +\frac{4\sigma^2}{n}\ln \left(\frac{2}{\beta}\right).
		\end{align*}
		Thus we have, under (\ref{12}),
		\begin{equation*}
		\mathbb{P}_{f}\left(\Phi_{K,\alpha} =0\right) \leq \beta.
		\end{equation*}
		Moreover, there exists some constant $\kappa>0$ such that, for every $K$ and $n \geq 32\ln (2/\alpha)$
		\begin{equation}\label{7}
		q^\alpha_{K,1-\beta/2} \leq \frac{2\kappa}{\sqrt{n(n-1)}} \ln \left(\frac{2}{\alpha}\right) \sqrt{\frac{2\int_{E^2}^{}K^2(x,y)d\nu(x)d\nu(y)}{\beta} }.
		\end{equation}
	\end{proposition}
	To prove the first part of this result, we simply use Markov's inequality for the term $T_K$ and an exponential inequality for non-central Chi-square variables due to (\cite{birge2001alternative}) for the term $\hat{\sigma}^2_n$. The control of $q^\alpha_{K,1-\beta/2}$ derives from a property of Gaussian chaoes combined with an exponential inequality (due to \cite{de2012decoupling} and \cite{huskova1993consistency}). The detailed proof is given in the Appendix.
	
	The following theorem gives a condition on $\|f\|^2$ for the test to be powerful.
	
	\begin{thm}\label{thmnormdiff}
		Let $\alpha,\ \beta$ be fixed levels in $(0,1)$, $\kappa$ be a positive constant, $K$ be a symmetric kernel function, and $\Phi_{K,\alpha}$ be the test defined by (\ref{testfunc}). Let $C_K$ be an upper bound for $\int_{E^2}^{}K^2(x,y)d\nu(x)d\nu(y)$. Then for all $n \geq 32\ln (2/\alpha)$, we have $\mathbb{P}_{f}(\Phi_{K,\alpha} = 0) \leq \beta$, as soon as
		\begin{align}\label{normformul}
			\| f\|^2  &\geq \|f-K[f] \|^2 + \frac{16\left( \norm{f}_\infty^2 +\sigma^2 \right)}{n\beta}\nonumber\\
			{} &+ \frac{4}{\sqrt{n(n-1)\beta}} \left(\kappa D_{n,\beta}\ln \left(\frac{2}{\alpha}\right) +\sqrt{2}\left(\norm{f}_\infty^2 +\sigma^2 \right) \right)\sqrt{C_K}.
		\end{align}
	\end{thm}
	The right hand side of the above inequality corresponds to a bias-variance trade-off. For particular choices of the kernel function $K$, these terms will be upper bounded in Section \ref{sect3}. 
	
	\subsection{Performance of the Monte Carlo approximation.}\label{monteSect}
	In this section, we introduce a Monte Carlo method used to approximate the conditional quantiles $q^{(X)}_{K,1-\alpha}$ by $\hat{q}^{(X)}_{K,1-\alpha}$ as follows. We consider the set of $2B$ independent sequences of i.i.d standard Gaussian variables
	\begin{equation*}
	\{\epsilon^b,\ 1\leq b \leq B \} \quad \text{and} \quad \{\epsilon^{'b},\ 1\leq b \leq B \},
	\end{equation*}
	where $\epsilon^b =\{\epsilon^b_i \}_{i=1}^n$, $\epsilon^{'b} =\{\epsilon^{'b}_i \}_{i=1}^n$, $1\leq b \leq B$.\\
	We define
	\begin{equation*}
	V_K^{\left(\epsilon^b, \epsilon^{'b} \right)} = \frac{\frac{1}{n(n-1)}\sum_{i\neq j=1}^{n}K(X_i,X_j)\epsilon_i^b \epsilon_j^b}{\frac{1}{n}\sum_{i=1}^{n/2}\left(\epsilon_{2i-1}^{'b} - \epsilon_{2i}^{'b} \right)^2 },
	\end{equation*}
	where $X=\left(X_1,\cdots,X_n\right)$ are observed from model (\ref{y1y2-equation}).\\
	Under $(H_0)$, conditionally on $X$, the variables $V_K^{\left(\epsilon^b, \epsilon^{'b} \right)}$ have the same distribution function as $V_K$ and as $V_K^{\left(0\right)}$. We denote by $F_{K,B}$ the empirical distribution function of the sample $\left\{V_K^{\left(\epsilon^b, \epsilon^{'b} \right)},\ 1\leq b \leq B \right\} $, conditionally on $X$.
	\begin{equation*}
	\forall x\in \mathbb{R},\quad F_{K,B}(x)=\frac{1}{B}\sum_{b=1}^{B}\mathbbm{1} \left\lbrace V_K^{\left(\epsilon^b, \epsilon^{'b} \right)} \leq x \right\rbrace .
	\end{equation*}
	Then the Monte Carlo approximation of $q^{(X)}_{K,1-\alpha}$ is defined by
	\begin{equation*}
	\hat{q}^{(X)}_{K,1-\alpha}=F^{-1}_{K,B}(\alpha)= \inf \left\lbrace t\in \mathbb{R},\ F_{K,B}(t) \geq 1-\alpha \right\rbrace.
	\end{equation*}
	We recall the test function defined in (\ref{testfunc}) and we reject $(H_0)$ when $V_K >q^{(X)}_{K,1-\alpha}$ with $q^{(X)}_{K,1-\alpha}$ the $(1-\alpha)$ quantile of $V_K^{\left(0\right)}$ defined by (\ref{test.epsilon}) conditionally on $X$. Now, by using the estimated quantile $\hat{q}^{(X)}_{K,1-\alpha}$, we consider the test given by
	\begin{equation}\label{test.mon}
	\widehat{\Phi}_{K,\alpha} = \mathbbm{1} \left\lbrace V_K > \hat{q}^{(X)}_{K,1-\alpha} \right\rbrace .
	\end{equation}
	For the test defined in (\ref{test.mon}), the probabilities of first and second kind errors can above upper bounded. This is the purpose of the two following propositions, whose proofs are given in \cite{fromont2013two}.
	
	\begin{proposition}\label{8}
		Let $\alpha$ be some fixed level in $(0,1)$, and $\widehat{\Phi}_{K,\alpha}$ be the test defined by (\ref{test.mon}). Then,
		\begin{equation*}
		\mathbb{P}_{(H_0)}\left(\widehat{\Phi}_{K,\alpha}=1 \bigg | X \right) \leq \frac{B\alpha +1}{B+1}.
		\end{equation*}
	\end{proposition}
	
	\begin{proposition}\label{9}
		Let $\alpha$ and $\beta$ be fixed levels in $(0,1)$ such that $\alpha_B = \alpha-\sqrt{\ln B/(2B)}$ and $\beta_B=\beta - 2/B>0$. Let $\widehat{\Phi}_{K,\alpha}$ be the test given in (\ref{test.mon}). Let $ A_K, B_K, D_{n,\beta}$ and $\kappa$ as in Proposition \ref{2nderror}, and let $q^{\alpha_B}_{K,1-\beta_B/2}$ be the $(1-\beta_B/2)$ quantile of $q^{(X)}_{K,1-\alpha_B}$. If
		\begin{equation}\label{10}
		\langle K[f],f \rangle > \sqrt{\frac{16A_K + 8B_K}{\beta}}+ D_{n,\beta_B}\  q^{\alpha_B}_{K,1-\beta_B/2},
		\end{equation}
		then $\mathbb{P}_{f}\left(\widehat{\Phi}_{K,\alpha}=0 \right)\leq \beta$. Moreover,
		\begin{equation}\label{11}
		q^{\alpha_B}_{K,1-\beta_B/2} \leq \frac{2\kappa}{\sqrt{n(n-1)}} \ln \left(\frac{2}{\alpha_B}\right) \sqrt{\frac{2\int_{E^2}^{}K^2(x,y)d\nu(x)d\nu(y)}{\beta_B} }.
		\end{equation}
	\end{proposition}
	\textit{Comments.} When comparing (\ref{10}) and (\ref{11}) with (\ref{12}) and (\ref{7}) in Proposition \ref{2nderror}, we notice that they asymptotically coincide when $B\rightarrow +\infty$. Moreover, if $\alpha=\beta=0.05$ and $B\geq 6000$, the multiplicative factor of $\kappa n\sqrt{B_K}$ is of order $1.2$ in (\ref{11}) compared with (\ref{7}).
	
	
	\section{Two particular examples of kernel function.}\label{sect3}
	In this section, we specify the performances of the above test for two examples of the kernels including projection kernels and Gaussian kernels. 
	\subsection{Projection kernels.}\label{sect.proj}
	We assume $E=[0,1]$. We consider the projection kernel defined in (\ref{projectionkernel}) and aim to give a more explicit formulation for the result of Theorem \ref{thmnormdiff} under the choice of this kernel. We also evaluate the uniform separation rates over Besov bodies.
	
	\begin{corollary}\label{col.f.gaussian}
		Let $\alpha,\ \beta \in (0,1)$ and $\kappa >0$ be a constant. Let $\Phi_{K,\alpha}$ be defined in (\ref{testfunc}), where $K$ is the projection kernel defined by (\ref{projectionkernel}). We denote by $S$ the linear subspace of $\mathbb{L}^2(\left[0,1\right],d\nu)$, generated by the functions $\{\phi_\lambda,\ \lambda \in \Lambda \}$, and we assume that the dimension of $S$ is equal $D$. Then $n \geq 32\ln (\alpha/2)$ if
		\begin{align*}
			\| f\|^2  &\geq \|f-\Pi_S(f) \|^2 + \frac{16\left(\norm{f}^2_\infty +\sigma^2 \right)}{n\beta}\\
		{} &+ \frac{4\sqrt{D}}{\sqrt{n(n-1)\beta}} \left(\kappa D_{n,\beta}\ln \left(\frac{2}{\alpha}\right) +\sqrt{2}\left(\norm{f}^2_\infty +\sigma^2 \right) \right),
		\end{align*}
		then
		\begin{equation*}
		\mathbb{P}_{f}\left(\Phi_{K,\alpha} = 0\right) \leq \beta.
		\end{equation*}
	\end{corollary}
	
	Let us consider the particular case when the kernel $K$ is the projection kernel onto the space generated by functions of the Haar basis defined as follows.\\
	Let $\{\phi_0,\ \phi_{\left(j,k\right)},\ j\in \mathbb{N}, k\in \{0,\cdots,2^j -1 \}$ be the Haar basis of $\mathbb{L}^2([0,1])$ with
	\begin{equation}
	\phi_0(x)=\mathbbm{1}_{[0,1]}(x) \quad \text{and} \quad \phi_{j,k}(x)=2^{j/2}\psi(2^j x-k),
	\end{equation}
	where $\psi(x)=\mathbbm{1}_{\left[0,1/2\right)}(x)-\mathbbm{1}_{\left[1/2,1\right]}(x)$. The linear subspace $S$ is generated by a subset of the Haar basis. More precisely, we denote by $S_0$ the subspace of $\mathbb{L}^2([0,1])$ generated by $\phi_0$, and we define 
	\begin{equation}\label{K.project.0}
	K_0(x,x^{'})=\phi_0(x)\phi_0(x^{'}).
	\end{equation}
	We also consider, for $J\geq 1$ the subspace $S_J$ generated by $\{\phi_\lambda,\ \lambda \in \{0\} \cup \Lambda_J\}$ with $\Lambda_J = \{(j,k),\ j\in \{0,\cdots,J-1 \},\ k\in \{0,\cdots,2^j-1 \} \}$, and
	\begin{equation}\label{K.project.J}
	K_J(x,x^{'})=\sum_{\lambda \in \{0\} \cup \Lambda_J}^{}\phi_\lambda(x) \phi_\lambda(x^{'}).
	\end{equation}
	We set $\alpha_0=\left\langle f,\phi_0 \right\rangle $ and for every $j\in \mathbb{N},\ k\in \{0,\cdots,2^j-1 \}$, $\alpha_{j,k}=\left\langle s,\phi_{j,k}\right\rangle $.\\
	We now introduce the Besov body defined for $\delta >0,\ R>0$ by
	\begin{align*}
	\mathcal{B}_{2,\infty}^\delta (R)=\left\lbrace f\in \mathbb{L}^2([0,1],d\nu),\ f=\alpha_0 \phi_0 +   \sum_{j\in \mathbb{N}}^{} \sum_{k=0}^{2^j-1} \alpha_{j,k}\phi_{j,k}\right. &/ \left. \alpha_0^2 \leq R^2, \ \forall j \in \mathbb{N}, \right. \\
	{} & \left. \sum_{k=0}^{2^j-1} \alpha^2_{j,k} \leq R^2 2^{-2j\delta} \right\rbrace .
	\end{align*}
	For all $J\geq 0$, we consider the kernel function $K_J$ defined by (\ref{K.project.0}), (\ref{K.project.J}) and the associated test function $\Phi_{K_J,\alpha}$ defined in (\ref{testfunc}) with $K=K_J$. For an optimal choice of $J$, realizing a good compromise between the bias term and the variance term appearing in (\ref{normformul}), we give a condition of $\|f\|^2$ for $f\in \mathcal{B}_{2,\infty}^\delta(R)$ which ensures that the power of our test is larger than $1-\beta$.

	\begin{proposition}\label{prop.project}
		Let $\alpha,\ \beta\in (0,1)$. For all $J\geq 0$, let $K_J$ defined by (\ref{K.project.0}), (\ref{K.project.J}) and consider the test function $\Phi_{{K_{J^{*}},\alpha}} =  \mathbbm{1} \{V_{K_{J^{*}}} > q^{(X)}_{K_{J^{*}},1-\alpha}\}$ where
		\begin{equation}\label{J.optimal}
		J^{*} = \left[\log_2\left(n^{2/(1+4\delta)} \right)\right].
		\end{equation}
		For all $f\in \mathcal{B}_{2,\infty}^\delta(R)$ such that
		\begin{equation}\label{normproj}
			\|f\|^2 \geq C(\alpha,\beta,\sigma,R,\|f\|_\infty) n^{-4\delta/(1+4\delta)},
		\end{equation}
		we have $\mathbb{P}_{f}(\Phi_{K_{J^{*}},\alpha} = 0) \leq \beta$. 
	\end{proposition}
	
	\textit{Comments.}
	\begin{enumerate}
		\item[1.] Non asymptotic lower bounds for the rates of testing in signal detection over Besov bodies are given in \cite{baraud2002non}. These lower bounds coincide with the bound given in (\ref{normproj}), hence our result is sharp.
		\item [2.] In (\ref{J.optimal}), $J^{*}$ depends on $\delta$, the regularility parameter of the Besov body, so it leads to the natural question of the choice if this parameter. In order to propose a procedure that is adaptive with respect to the regularity of the unknown regression function $f$, we introduce aggregated tests in Section \ref{Sec4}.
	\end{enumerate}

	\subsection{Gaussian kernels.}\label{sect.gauss}
	For this second example, we assume that $E=\mathbb{R}$. We consider the Gaussian kernel defined in (\ref{13}) and rewrite the result of Theorem \ref{thmnormdiff} under the choice of this kernel. We also evaluate the uniform separation rates over Sobolev balls for this test. 
	\begin{corollary}\label{col.gauss}
		Let $\alpha, \beta \in (0,1)$, $\kappa>0$ be a constant and $\Phi_{K,\alpha}$ be the test function defined in (\ref{testfunc}) where $K$ is defined in (\ref{13}). For $n \geq 32\ln (2/\alpha)$ if
		\begin{align}\label{normgaus}
			\| f\|^2  &\geq \|f-k_h*f \|^2 + \frac{16\left(\norm{f}^2_\infty +\sigma^2 \right)}{n\beta} \nonumber \\
			{} &+ \frac{4\norm{\nu}_\infty}{(2\pi)^{1/4}\sqrt{n(n-1)\beta h}} \left(\kappa D_{n,\beta}\ln \left(\frac{2}{\alpha}\right) +\sqrt{2}\left(\norm{f}^2_\infty +\sigma^2 \right) \right).
		\end{align}
		We obtain that
		\begin{equation*}
		\mathbb{P}_{f}(\Phi_{K,\alpha} = 0) \leq \beta.
		\end{equation*}
	\end{corollary}
	Let $E=\mathbb{R}$ and $\mathcal{L}=\mathbb{N}^{*}$. For $x,\ y$ in $\mathbb{R}$ and $h=2^{-l}$, for all $l\in \mathcal{L}$, we consider
	\begin{equation}\label{K.m}
	K_l(x,y) =\frac{1}{2^{-l}} k\left(\frac{x-y}{2^{-l}} \right),
	\end{equation}
	with 
	\begin{equation*}
	k(u) = \frac{1}{\sqrt{2\pi}}\exp{\left(-\frac{u^2}{2} \right)}.
	\end{equation*}
	Let us introduce for $\delta>0$ the Sobolev ball $\mathcal{S}^\delta(R)$ defined by 
	\begin{equation*}
	\mathcal{S}_\delta(R) = \left\{s: \mathbb{R} \rightarrow \mathbb{R}\ \big/ \ s\in \mathbb{L}^1(\mathbb{R})\cap \mathbb{L}^2(\mathbb{R}),\ \int_{\mathbb{R}}^{} |u|^{2\delta} |\hat{s}(u)|^2 du \leq 2\pi R^2 \right\},
	\end{equation*}
	where $\hat{s}$ denotes the Fourier transform of $s$: $\hat{s}(u)=\int_{\mathbb{R}}^{}s(x)e^{i\langle x,u \rangle}dx$.\\
	For all $l\in \mathcal{L}$, we consider the kernel function $K_l$ defined by (\ref{K.m}) and the associated test function $\Phi_{K_l,\alpha}$ defined in (\ref{testfunc}) with $K=K_l$. For an optimal choice of $l$, realizing a good compromise between the bias term and the variance term appearing in (\ref{normgaus}), we give a condition of $\|f\|^2$ for $f\in \mathcal{S}_\delta(R)$ which ensures that the power of our test is larger than $1-\beta$.
	
	\begin{proposition}\label{prop.gauss}
		Let $\alpha,\ \beta \in (0,1)$. For all $l\in \mathcal{L}$, let $K_l$ defined by (\ref{K.m}) and the test function $\Phi_{K_l,\alpha} = \mathbbm{1} \{V_{K_l} > q^{(X)}_{K_l,1-\alpha}\}$ we set 
		\begin{equation}\label{m.optimal}
		l^{*} = \left[\log_2\left(n^{\frac{2}{1+4\delta}} \right) \right].
		\end{equation}
		For all $f\in \mathcal{S}_\delta(R)$ such that
		\begin{equation}\label{funct26}
		\|f\|^2 \geq C(\alpha,\beta,\sigma,R,\|f\|_\infty)n^{-4\delta/(1+4\delta)}.
		\end{equation}
		We have $\mathbb{P}_{f}(\Phi_{K_{l^{*}},\alpha} = 0) \leq \beta$.
	\end{proposition}
	
	\textit{Comments.}
	\begin{enumerate}
		\item[1.] As in Proposition \ref{prop.project}, we obtain in the right hand term of (\ref{funct26}) a classical bound for the separation rates of testing over regular classes of alternatives such as Holderian balls (see \cite{ingster1993asymptotically}) for nonparametric minimax rates of testing in various setups.
		\item[2.] Non asymptotic lower bounds for the rates of testing in signal detection over Sobolev balls are given in \cite{fromont2006adaptive}. These bounds coincide with the bound given in (\ref{funct26}).
		\item [3.] In (\ref{m.optimal}), as previously, $l^{*}$ depends on $\delta$, the regularity parameter of the Sobolev ball, so it leads to the natural question of the choice of this parameter answered through the aggregated tests in Section \ref{Sec4}.
	\end{enumerate}
	
	\section{Multiple or aggregated tests based on collections of kernel \label{Sec4} functions.}\label{SectAggregateTest}
	In the previous section, we have considered testing procedures based on a single kernel function $K$. However, the following question is natural: how can we choose the kernel, and its parameters. For instance, the orthonormal family in the projection kernel in Section \ref{sect.proj}, the bandwidth $h$ in the Gaussian kernel in Section \ref{sect.gauss}. \cite{baraud2003adaptive} proposed adaptive testing procedures based on the aggregation of a collection of tests. This idea is presented in a series of papers, among which \cite{fromont2013two} proposed an aggregation procedure. Following this idea, we consider in this section a collection of kernel functions instead of a single one. Beside that, we define a multiple testing procedure by aggregating the corresponding single tests, with an adapted choice of the critical values.
	
	\subsection{The aggregated testing procedure.}\label{Sect4.1}
	Let us describe the aggregated testing procedure by introducing a finite collection $\left\lbrace K_m,\ m\in \mathcal{M} \right\rbrace $ of symmetric kernel functions: $E\times E \rightarrow \mathbb{R}$. For  $m\in \mathcal{M}$, we replace $K$ in (\ref{Vk}) and (\ref{test.epsilon}) by $K_m$ to define $V_{K_m}$ and $V_{K_m}^{\left(0\right)}$ and let $\left\lbrace w_m,\ m\in \mathcal{M} \right\rbrace $ be a collection of positive numbers such that $\sum_{m\in \mathcal{M}}^{} e^{-w_m} \leq 1$. Conditionally on $X$, for $u\in (0,1)$, we denote by $q^{(X)}_{m,1-u}$ the $(1-u)$ quantile of $V_{K_m}^{\left(0\right)}$. Given $\alpha$ in $(0,1)$, we consider the test which rejects $(H_0)$ when there exists at least one $m$ in $\mathcal{M}$ such that
	\begin{equation*}
	V_{K_m} > q^{(X)}_{m,1-u^{(X)}_\alpha e^{-w_m}},
	\end{equation*}
	where $u_\alpha^{(X)}$ is defined by
	\begin{equation}\label{ualphaX}
	u_\alpha^{(X)} = \sup \left\lbrace u>0,\ \mathbb{P}\left(\sup_{m\in \mathcal{M}}\left(V_{K_m}-q^{(X)}_{m,1-u e^{-w_m}} \right) >0 \bigg |X \right)\leq \alpha \right\rbrace .
	\end{equation}
	We consider the test function $\Phi_{\alpha}$ defined by
	\begin{equation}\label{Phi.multi}
	\Phi_{\alpha} = \mathbbm{1}\left\lbrace \sup_{m\in \mathcal{M}}\left(V_{K_m}-q^{(X)}_{m,1-u_\alpha^{(X)} e^{-w_m}} \right)>0 \right\rbrace .
	\end{equation}
	Using the Monter Carlo method, we can estimate $u_\alpha^{(X)}$ and the quantiles $q^{(X)}_{m,1-u_\alpha^{(X)} e^{-w_m}}$ for all $m\in \mathcal{M}$. The following theorem provides a coltrol of the first and second kind error for the test $\Phi_\alpha$. The detailed proof is given in the Appendix.
	\begin{thm}\label{them4.1them}
		Let $\alpha,\beta$ be fixed levels in $(0,1)$ and $\Phi_{\alpha}$ be the test defined by (\ref{Phi.multi}). We have
		\begin{equation}
		\mathbb{P}_{(H_0)}\left(\Phi_{\alpha}=1 \right) \leq \alpha.
		\end{equation}
		And for all regression function $f$, we have
			\begin{equation}
		\mathbb{P}_{f}\left(\Phi_{\alpha}=0 \right) \leq \beta,
		\end{equation}
		as soon as there exists $m$ in $\mathcal{M}$ such that
		\begin{equation*}
		\mathbb{P}\left(V_{K_m} \leq q^{(X)}_{K_m,1-\alpha e^{-w_m}}\right) \leq \beta.
		\end{equation*}
	\end{thm}
	\textit{Comments}. This theorem shows that the aggregated test is of level $\alpha$, for all $n$. Moreover, as soon as the second kind error is controlled by $\beta$ for at least one test in the collection, the same holds for the aggregated procedure with the price that the level $\alpha$ is replaced by $u_\alpha^{(X)}e^{-w_m}$ to guarantee that the aggregated procedure is of level $\alpha$.
	\subsection{The aggregation of projection kernels.}\label{aggre.proj}
	Let us specify the performance of the aggregated test for a collection of projection kernels.
	\begin{corollary}\label{colaggreProj}
		Let $\alpha, \beta$ be fixed levels in $(0,1)$. Let $\left\lbrace S_m, m\in \mathcal{M} \right\rbrace $ be a finite collection of linear subspaces of $\mathbb{L}^2([0,1],d\nu )$, generated by the functions $\left\lbrace \phi_\lambda, \lambda \in \Lambda_m \right\rbrace $ and we assume that the dimension of $S_m$ is equal to $D_m$. We set, for all $m\in \mathcal{M}$, $K_m(x,y)=\sum_{\lambda \in \Lambda_m}^{}\phi_\lambda (x)\phi_\lambda (y)$. Let $\Phi_\alpha$ be defined by (\ref{Phi.multi}) with the collection of kernels $\left\lbrace K_m,\ m\in \mathcal{M} \right\rbrace $ and the collection $\left\lbrace w_m,m\in \mathcal{M} \right\rbrace $ of positive numbers such that $\sum_{m\in \mathcal{M}}^{}e^{-w_m} \leq 1$.
		
		Then $\Phi_\alpha$ is a level $\alpha$ test. Moreover, $\mathbb{P}_{f} \left(\Phi_\alpha =0 \right) \leq \beta$ if
		\begin{align}\label{aggreproj}
		\| f\|^2  &\geq \inf_{m\in \mathcal{M}} \left\lbrace  \|f-\Pi_{S_m}(f) \|^2 + \frac{16\left(\norm{f}^2_\infty +\sigma^2 \right)}{n\beta} \right. \nonumber \\
		{} &+ \left.  \frac{4\sqrt{D_m} }{\sqrt{n(n-1)\beta}} \left(\kappa D_{n,\beta} \left( \ln \left(\frac{2}{\alpha}\right)+w_m\right)  +\sqrt{2}\left(\norm{f}^2_\infty +\sigma^2 \right) \right)\right\rbrace ,
		\end{align}
		where $\kappa >0$ and $n \geq 32\ln (\alpha/2)$.
	\end{corollary}

	\textit{Comments}. Comparing this result with the one obtained in Corollary \ref{col.f.gaussian} for the single test based on a projection kernel, we can see that the multiple testing procedure allows to obtain the infimum over all $m$ in $\mathcal{M}$ in the right hand side of (\ref{aggreproj}) at the price of the additional term $w_m$. 
	
	Let us consider the particular case when the collection of kernels $\left\lbrace K_m,\ m\in \mathcal{M} \right\rbrace $ is the collection of projection kernels based on the constructions in Section \ref{sect.proj}. Let for some $\bar{J} \geq 1,\ \mathcal{M}_{\bar{J}} =\left\lbrace J,\ 0\leq J\leq \bar{J} \right\rbrace $, and for all $J$ in $\mathcal{M}_{\bar{J}}$, $w_J = 2\left(\ln(J+1)+\ln (\pi/\sqrt{6}) \right)$. 
	
	We consider $\Phi_\alpha^{(1)}$, the test defined by (\ref{Phi.multi}) with the collection of kernels $\left\lbrace K_J, 0\leq J\leq \bar{J} \right\rbrace $ where $K_0,\ K_J, 0<J\leq \bar{J}$ defined in (\ref{K.project.0}), (\ref{K.project.J}). We obtain from the Corollary \ref{colaggreProj} that there exists some constant $C(\alpha,\beta,\sigma,\|f\|_\infty) $ such that $\mathbb{P}_{f} \left(\Phi_\alpha^{(1)} = 0 \right) \leq \beta$ as soon as
	\begin{equation}\label{aggrePrker}
	\| f\|^2  \geq C(\alpha,\beta,\sigma,\|f\|_\infty) \inf_{J\in \mathcal{M}_{\bar{J}}} \left\lbrace  \|f-\Pi_{S_J}(f) \|^2 +\ln(J+2)\frac{2^{J/2}}{n} \right\rbrace .
	\end{equation}
	For any $\delta >0, R, R^{'} >0$ we consider 
	\begin{equation}
	\mathcal{B}_{2,\infty}^\delta(R,R^{'}) = \left\lbrace f \ : f\in \mathcal{B}_{2,\infty}^\delta(R),\  \norm{f}_\infty \leq R^{'}  \right\rbrace .
	\end{equation}
	\begin{corollary}\label{corol4.3}
		Let $\alpha, \beta \in (0,1)$. For all $J \in \mathcal{M}_{\bar{J}}$, we consider the test function $\Phi_\alpha^{(1)}$. Assuming that $\ln \ln(n) \geq 1,\ 2^{\bar{J}}\geq n^2$. Then, for any $\delta, R, R^{'} >0 $ we set
		\begin{equation*}
		J^{**} =  \left[ \log_2\left(\left(\frac{n}{\ln \ln(n)} \right)^{\frac{2}{4\delta+1}} \right) \right] .
		\end{equation*}
		 For all  $f\in \mathcal{B}_{2,\infty}^\delta(R,R^{'}) $ such that
		 \begin{equation}\label{col4.3}
		 \norm{f}^2 \geq C(\alpha,\beta,\sigma,R,R^{'})\left(\frac{\ln \ln(n)}{n}\right)^{\frac{4\delta}{4\delta+1}},
		 \end{equation}
		 we have $\mathbb{P}_{f}\left(\Phi_\alpha^{(1)}=0 \right) \leq \beta$.
	\end{corollary}

	\textit{Comments.} We obtain a right hand term in (\ref{col4.3}) of order $\left(\ln\ln(n)/n \right)^{4\delta/(1+4\delta)}$. This rate of testing was shown to be optimal for the signal detection in a Gaussian white noise by \cite{spokoiny1996adaptive}. In particular, he showed that the logarithm factor is the price to pay for adaptation.
	
	\subsection{The aggregation of Gaussian kernels.}\label{aggre.gauss}
	We here consider the aggregated test based on a collection of Gaussian kernels.
	\begin{corollary}\label{gauss.aggre}
		Let $\alpha, \beta \in (0,1)$, $\left\lbrace h_l,\ l\in \mathcal{L} \right\rbrace $ be a collection of positive bandwidths, we consider $\left\lbrace K_l, l\in \mathcal{L} \right\rbrace $ a collection of Gaussian kernels corresponding to the above collection of positive bandwidths, where $K_l$ defined in (\ref{K.m}). Let $\Phi_\alpha$ be defined by (\ref{Phi.multi}) with the collection of kernel $\left\lbrace K_l,\ l\in \mathcal{L} \right\rbrace$ and a collection $\left\lbrace w_l,\ l\in \mathcal{L} \right\rbrace $ of positive numbers such that $\sum_{l\in \mathcal{L}}^{}e^{-w_l} \leq 1 $.\\
		Then $\Phi_\alpha$ is a level $\alpha$ test. Moreover, there exists $\kappa>0$ such that if
		\begin{align}
		\| f\|^2  &\geq \inf_{l\in \mathcal{L}} \left\lbrace  \|f-k_l*f \|^2 + \frac{16\left(\norm{f}^2_\infty +\sigma^2 \right)}{n\beta} \right. \nonumber\\
		{} &+ \left.  \frac{4\norm{\nu}_\infty}{(2\pi)^{1/4}\sqrt{n(n-1)\beta h_l}} \left(\kappa D_{n,\beta} \left( \ln \left(\frac{2}{\alpha}\right) +w_l \right) +\sqrt{2}\left(\norm{f}^2_\infty +\sigma^2 \right) \right) \right\rbrace ,
		\end{align}
		We obtain that $\mathbb{P}_{f}(\Phi_{K,\alpha} = 0) \leq \beta.$
	\end{corollary}
	For $l\in \mathcal{L}=\mathbb{N}\setminus \{0\}$. We consider the particular case where we take $h_l=2^{-l}$ and $w_l=2\left(\ln(l+1)+\ln(\pi^2/6) \right)$ for all $l\in \mathcal{L}$. Let $\Phi_\alpha^{(2)}$ be the test defined by (\ref{Phi.multi}) with the collection of Gaussian kernels $\left\lbrace K_l, l\in \mathcal{L} \right\rbrace$ and $\left\lbrace w_l, l\in \mathcal{L} \right\rbrace $. We obtain from Corollary \ref{gauss.aggre} that there exists $C(\alpha,\beta,\sigma,\norm{f}_\infty)$ such that $\mathbb{P}_{f}\left(\Phi_\alpha^{(2)} =0 \right) \leq \beta$ if
	\begin{equation}\label{gaussiAggre}
	\norm{f}^2 \geq C(\alpha,\beta,\sigma,\norm{f}_\infty) \inf_{l\in \mathcal{L}} \left\lbrace  \|f-k_l*f \|^2 + \frac{w_l}{n\sqrt{2^{-l}}} \right\rbrace .
	\end{equation}
	For $\delta >0, R, R^{'} >0$ we consider
	\begin{equation}
	\mathcal{S}_\delta(R,R^{'}) = \left\lbrace f: f\in \mathcal{S}_\delta(R), {\norm{f}_\infty} \leq R^{'} \right\rbrace. 
	\end{equation}
	\begin{corollary}
		Let $\alpha, \beta \in (0,1)$. For all $l\in \mathcal{L}$, we consider the test function $\Phi_\alpha^{(2)}$ and assume that $\ln \ln \geq 1$. For any $\delta >0, R, R^{'} >0$, we set
		\begin{equation*}
		l^{**} =\left[  \log_2\left(\left(\frac{n}{\ln \ln(n)} \right)^{\frac{2}{4\delta+1}} \right) \right]  .
		\end{equation*}
	For all $f\in \mathcal{S}_\delta(R,R^{'})$ such that 
	\begin{equation*}
	\norm{f}^2 \geq C(\alpha,\beta,\sigma,R,R^{'})\left(\frac{\ln \ln(n)}{n}\right)^{\frac{4\delta}{4\delta+1}},
	\end{equation*}
	 we have $\mathbb{P}_{f}\left(\Phi_\alpha^{(2)}=0 \right) \leq \beta$.
	\end{corollary}

\textit{Comments.} The rate of testing is of order $\left(\ln \ln(n)/n \right)^{4\delta/(1+4\delta)}$. This rate was shown to be optimal over periodic Sobolev balls up to the logarithm, by \cite{castillo2006exact}.
	\section{Simulation study.}\label{simulationSect}
	\subsection{Presentation of the simulation study.}
	 We study our aggregated testing procedures from a practical point of view in this section. We consider $E=[0,1],\ n=100$ and choose $\alpha=0.05$. In the following simulation, $X_1,\cdots,X_n$ are i.i.d uniform random variables on $[0,1]$.
	
	Let us introduce the collection of symmetric kernel functions and the aggregated testing procedure $\Phi_\alpha$ defined by (\ref{Phi.multi}) as follows. First, we consider the test $\Phi_\alpha^{(1)}$ denoted by P corresponding to a collection of projection kernels. To be more explicit, we consider the Haar basis $\{\phi_0,\ \phi_{\left(j,k\right)},\ j\in \mathbb{N}, k\in \{0,\cdots,2^j -1 \}$ introduced in Section \ref{sect.proj}. Let $K_0(x,x^{'})=\phi_0(x)\phi_0(x^{'})$ and for $J\geq 1$ $K_J(x,x^{'})=\sum_{\lambda \in \{0\} \cup \Lambda_J}^{}\phi_\lambda(x) \phi_\lambda(x^{'})$ with $\Lambda_J = \{(j,k),\ j\in \{0,\cdots,J-1 \},\ k\in \{0,\cdots,2^j-1 \} \}$. Let $\mathcal{M}_{\overline{J}}=\left\lbrace J,\ 0\leq J\leq 7 \right\rbrace $ and for all $J$ in $\mathcal{M}_{\overline{J}}$, $w_J = 2\left(\ln(J+1)+\ln (\pi/\sqrt{6}) \right)$. We consider $\Phi_\alpha^{(1)}$ the multiple testing procedure with the collection of kernels $\left\lbrace K_J,\ J\in \mathcal{M}_{\overline{J}} \right\rbrace $. 
	
	Second, we also consider the multiple test associated with the collection of Gaussian kernel functions defined in Section \ref{aggre.gauss}. For $\mathcal{L}=\left\lbrace 1,2,\cdots,6 \right\rbrace$ we take $\left\lbrace h_l, l\in \mathcal{L} \right\rbrace = \left\lbrace 1/24, 1/16, 1/12, 1/8, 1/4, 1/2 \right\rbrace $, let $K_l(x,y) = \frac{1}{h_l}k\left(\frac{x-y}{h_l} \right)$ with $k(u)=(2\pi)^{-1/2}\exp \left(-u^2/2 \right)$. Then taking $w_l=1/|\mathcal{L}|=1/6$, we consider $\Phi_\alpha^{(2)}$ the multiple testing procedure denoted by G, with the collection of kernels $\left\lbrace K_l, l\in \mathcal{L} \right\rbrace $. 
	
	At last, we are interested in the collection of both projection and Gaussian kernels. We define $\Phi_\alpha^{(3)}$ denoted by PG, the multiple testing procedure with the collection of kernels $\left\lbrace  K_p,\ p\in \mathcal{P}=\mathcal{M}_{\overline{J}} \cup \mathcal{L}  \right\rbrace  $. For $p\in  \mathcal{M}_{\overline{J}}$ we take $w_{p}=\ln(J+1)+\ln (\pi/\sqrt{6})$ and for $p\in \mathcal{L}$ we take $w_{p}=1/12$.
	
	We recall that the test  rejects $(H_0)$ when there exists at least one $m$ in $\mathcal{M}$ such that $V_{K_m} > q^{(X)}_{m,1-u^{(X)}_\alpha e^{-w_m}}$. Hence, for each observation $X=\left(X_1,\cdots,X_n \right)$ we have to estimate $u_\alpha^{(X)}$ defined by (\ref{ualphaX}) and $q^{(X)}_{m,1-u^{(X)}_\alpha e^{-w_m}}$. Applying the Monte Carlo method introduced in the Section \ref{monteSect}, these quantities are well approximated. To be more explicit, we generate $400000$ samples of $\left\lbrace \epsilon^b \right\rbrace_{b=1}^{400000} $ and $\left\lbrace \epsilon^{'b} \right\rbrace_{b^{'}=1}^{400000} $, in which we use one half to approximate the conditional probability occurring in (\ref{ualphaX}) and other half is used to estimate the distribution of each $V_{K_m}^{\left(0\right)}$. We note that $u_\alpha^{(X)}$ is approximated by taking $u$ in a regular grid of $[0,1]$ with bandwidth $2^{-16}$ and choosing the approximation of $u_\alpha^{(X)}$ as the largest value of the grid such that the estimated conditional probabilities in (\ref{ualphaX}) are less than $\alpha$.
	\subsection{Simulation results.}
	We first study the probabilities of first kind error of each test. We realize $5000$ simulations of $X$. For each simulation, we determine the conclusions of tests P, G and PG where the critical values are approximated by the Monte Carlo methods described above. The probabilities of first kind error of tests are estimated by the number of rejections for these tests divided by $5000$. The obtained estimated levels of tests and the corresponding confidence intervals (CI) are showed in the Table \ref{tab1}.
	
	\begin{table}
		\begin{center}
			\begin{tabular}{|l|c|c|}
			\hline
			& \textbf{1st error} & \textbf{CI}\\
			\hline
			P & 0.0504&$\left[0.033,0.068\right]$ \\
			\hline
			G & 0.0506 &$\left[0.032,0.068\right]$ \\
			\hline
			PG &0.0498 &$\left[0.032,0.0657\right]$ \\
			\hline
			\end{tabular}
		\caption{\label{tab1}The probabilities of first kind error of the test for $\alpha=0.05$ and the upper and lower bounds of an asymptotic confidence interval with confidence level $99\%$.}
		\end{center}
	\end{table}
	
	We then study the probabilities of rejection for each test under several alternatives. We first consider the following alternative,
	\begin{equation*}
	f_{1,a,\epsilon}(x)=\epsilon\mathbbm{1}_{\left[0,a \right)}(x)-\epsilon\mathbbm{1}_{\left[a,2a \right)}(x)  ,
	\end{equation*}
	with $0 <\epsilon \leq 1 $ and $0<a<1$.
	Second, we consider the alternative defined by
	\begin{equation*}
	f_{2,\tau}(x) = \tau \sum_{j}^{} \frac{h_j}{2} \left(1 + \text{sgn}\left(x-p_j \right) \right),
	\end{equation*}
	with $\tau >0$,  and $h_j \in \mathcal{Z},\ 0 <p_j<1$ for all $j$. Next, we consider the following alternative,
	\begin{equation*}
	f_{3,c}(x) = c \cos (10\pi x),
	\end{equation*}
	with $c>0$. The last alternative, for which we aim to compare our results with the results of \cite{eubank1993testing} is defined as follows
		\begin{equation*}
	f_{4,\varrho,j}(x) = \varrho \cos(2\pi j x),
	\end{equation*}
	where $\varrho \geq 0$ and $j\in \mathbb{N}^{*}$. 
	
		\begin{table}[!htb]
		\begin{center}
			\scalebox{0.9}{
			\begin{tabular}{|l|c|c|c|c|c|c|c|c|}
				\hline
				\multirow{2}{*}{$\left(a,\epsilon\right)$}& \multicolumn{2}{|c|}{$\left(1/4,0.7\right)$}& \multicolumn{2}{|c|}{$\left(1/4,0.9\right)$}
				& \multicolumn{2}{|c|}{$\left(1/4,1\right)$}& \multicolumn{2}{|c|}{$\left(1/8,1\right)$}\\
				\cline{2-9}
				
				&$\hat{p}$&CI&$\hat{p}$&CI&$\hat{p}$&CI&$\hat{p}$&CI\\
				\hline
				P &0.876&$\left[0.849,0.903 \right] $&0.986&$\left[0.976,0.996 \right] $&0.996& $\left[0.990 ,1.001 \right]$&0.699&$\left[0.662,0.736\right]$\\
				\hline
				
				G &0.831& $\left[0.801,0.861\right]$&0.977&$\left[0.965,0.989\right]$&0.992&$\left[0.985,0.999\right]$&0.635&$\left[0.596,0.674\right]$\\
				\hline
				
				PG &0.884&$\left[0.858,0.910\right]$&0.984&$\left[0.973,0.994\right]$&0.996&$\left[0.991,1.001\right]$&0.690&$\left[0.652,0.727\right]$\\
				\hline
				
			\end{tabular}}
			\caption{\label{tab2}The power of the test for the alternative $f_{1,a,\epsilon}$ corresponding to $\left(a,\epsilon\right)=\left(1/4,0.7\right),\ \left(1/4,0.9\right),\ \left(1/4,1\right), \left(1/8,1\right)$ and the upper and lower bounds of an asymptotic confidence intervals with confidence level $99\%$.}
		\end{center}
	\end{table}
	
	\begin{table}[!htb]
		\begin{center}
			\begin{tabular}{|l|c|c|c|c|c|c|}
				\hline
				\multirow{2}{*}{$\tau$}& \multicolumn{2}{|c|}{$0.05$}& \multicolumn{2}{|c|}{$0.1$}
				& \multicolumn{2}{|c|}{$0.5$}\\
				\cline{2-7}
				&$\hat{p}$&CI&$\hat{p}$&CI&$\hat{p}$&CI\\
				\hline
				P &0.218&$ \left[0.177,0.243\right]$&0.654& $\left[0.615,0.693\right]$&1&* \\
				\hline
				
				G &0.208&$ \left[0.175,0.241\right]$&0.668&$\left[0.629,0.704\right]$&1&*\\
				\hline
				
				PG &0.210&$\left[0.177,0.243\right]$&0.678&$\left[0.639,0.716\right]$&1&
				*\\
				\hline
				
			\end{tabular}
			\caption{\label{tab3}The power of the test for the alternative $f_{2,\tau}$ corresponding to $\tau=1, 2, 3$ and the upper and lower bounds of an asymptotic confidence intervals with confidence level $99\%$.}
		\end{center}
	\end{table}
	
	\begin{table}[!htb]
		\begin{center}
			\begin{tabular}{|l|c|c|c|c|c|c|}
				\hline
				\multirow{2}{*}{$c$}& \multicolumn{2}{|c|}{$1$}& \multicolumn{2}{|c|}{$2$}
				& \multicolumn{2}{|c|}{$3$}\\
				\cline{2-7}
				&$\hat{p}$&ICI&$\hat{p}$&CI&$\hat{p}$&CI\\
				\hline
				P &0.35&$\left[0.311,0.389\right]$ &0.90&$\left[0.876,0.924\right] $&0.98&$\left[0.969,0.991\right] $\\
				\hline
				
				G &0.56&$\left[0.519, 0.600\right] $&0.98&$\left[0.967,0.991\right]$&1&*\\
				\hline
				
				PG &0.34&$\left[0.301,0.379\right]$&0.89&$\left[0.864,0.915\right]$&1&*\\
				\hline
				
			\end{tabular}
			\caption{\label{tab4}The power of the test for the alternative $f_{3,c}$ corresponding to $c=1,2,3$ and the upper and lower bounds of an asymptotic confidence intervals with confidence level $99\%$.}
		\end{center}
	\end{table}

	\begin{table}[!h]
		\begin{center}
			\scalebox{1}{
			\begin{tabular}{|l|c|c|c|c|c|}
				\hline
				&\text{Test}&$\varrho=0$&$\varrho=0.5$
				&$\varrho=1$& $\varrho=1.5$\\
				\hline 
				
				\multirow{3}{*}{}& P&0.049& 0.606 & 1 &1 \\
				$j=1$& G & 0.048 & 0.459 &0.99 &1\\
				& PG &0.048&0.441&0.99&1\\
				&EL1&0.074&0.837&1&1\\
				&EL2&0.062&0.805&1&1\\
				\hline
				\multirow{3}{*}{}& P&0.053& 0.224 &0.905&1\\
				$j=3$& G & 0.053 &  0.630 &0.922 &1 \\
				& PG &0.049&0.228& 1 &1\\
				&EL1&0.069&0.718&1&1\\
				&EL2&0.058&0.693&1&1\\
				\hline
				\multirow{3}{*}{}& P&0.043&  0.134& 0.696&0.990\\
				$j=6$& G & 0.044&  0.146 & 0.741 &0.995\\
				& PG &0.045&0.134&0.700&0.996 \\
				&EL1&0.076&0.134&0.428&0.979\\
				&EL2&0.056&0.107&0.368&0.961\\
				\hline
			\end{tabular}}
			\caption{\label{tab5}The power of the test for the alternative $f_{4,\varrho,j}$ corresponding to $\varrho=0,0.5,1,1.5,\ j=1,2,3$.}
		\end{center}
	\end{table}
	For each alternative $f$, we realize 1000 simulations of $X$. For each simulation, we determine conclusions of tests P, G and PG, where the critical values of our tests are still approximated by the Monte Carlo method. The powers of tests are estimated by the number of rejections divided by 1000. The obtained estimated powers of tests and lower bounds of the asymptotic confidence intervals with the confidence level $99 \%$ are represented in the Table \ref{tab2}, \ref{tab3} and \ref{tab4}. Table \ref{tab5} is proposed for comparing our tests and the two of tests $T_{nm}$ denoted by EL1,  $T_{n\lambda}$ denoted by EL2, which were proposed in \cite{eubank1993testing}. We recall briefly tests $T_{nm}$, $T_{n\lambda}$ as follows.
	\begin{equation*}
	T_{nm}=\frac{\sum_{|j|\leq m}^{'}|\tilde{a}_{jn}|^2-2m\sigma^2 }{2\sigma^2\sqrt{m}},
	\end{equation*}
	and
	\begin{equation*}
	T_{n\lambda}=\frac{n\sum_{|j|\leq (n-1)/2}^{'}|a_{jn}|^2\left(1+\lambda(2\pi j)^4 \right)^{-2} - 2\sigma^2 \sum_{|j|\leq (n-1)/2}^{'}(1+\lambda (2\pi j)^4 )^{-2}  }{2\sigma^2 \left(\sum_{|j|\leq (n-1)/2 }^{'} \left(1+\lambda(2\pi j)^4 \right)^{-4} \right)^{1/2} },
	\end{equation*}
	where $\sum_{}^{'}$ indicates summation excluding the zero index and $\tilde{a}_{jn}$ are the sample Fourier coefficients,
	\begin{equation*}
	\tilde{a}_{jn}=\frac{1}{n}\sum_{i=1}^{n}y_{i}e^{-2\pi j(i-1)/n}.
	\end{equation*}
	
	In the three alternatives $f_{1,a,\epsilon}$, $f_{2,\tau}$ and $f_{3,c}$, the test PG is more powerful than P and G tests. Our conclusion is that the test PG is a good choice in practice. In Table \ref{tab5}, we see in the firt column ($\rho=0$), which corresponds to the null hypothesis that our test is of level $\alpha=0.05$, which is not the case for the tests proposed by \cite{eubank1993testing}, which are only asymptotically of level $\alpha$. This explain why our test is generally less powerful than the tests EL1 and EL2 for $\rho=0.5$. In the other cases, we obtain as good or better results.

	\begin{appendices}
		
		\section{Proof of Proposition \ref{2nderror} }
		Let us prove the first part of Prop \ref{2nderror}. Recall that $q^\alpha_{K,1-\beta/2}$ denotes the $(1-\beta/2)$ quantile of $q^{(X)}_{K,1-\alpha}$ which is the $(1-\alpha)$ quantile of $V_K^{\left(0\right)}$ conditionally on $X$. We here want to find a condition on $\varepsilon_K = \mathbb{E}(T_K)$, ensuring that
		\begin{equation*}
		\mathbb{P}_{f}\left(V_K \leq q^\alpha_{K,1-\beta/2} \right) \leq \beta/2.
		\end{equation*}
		From Markov's inequality, we have for all $\lambda>0$
		\begin{equation} \label{P.H1}
			\mathbb{P}_{f} \left(|-T_K + \varepsilon_K| \geq \lambda \right) \leq \frac{\mathbb{E}(T_K^2) - \varepsilon_K^2}{\lambda^2}.
		\end{equation}
		Let us compute $\mathbb{E}\left(T_K^2 | X\right)$. We see that
		\begin{align*}
			T_K^2 &= \frac{1}{n^2(n-1)^2}\left[\sum_{\substack{i\neq j\neq u\neq v\\i,j,u,v=1}}^{n} K_{ij}K_{uv} Y_i Y_j Y_u Y_v + 4\sum_{\substack{i\neq j\neq u\\i,j,u=1}}^{n} K_{ij}K_{iu}Y_i^2Y_j Y_u +2\sum_{\substack{i\neq j\\i,j=1}}^{n} K^2_{ij} Y_i^2 Y_j^2 \right].
		\end{align*}
		Then
		\begin{align*}
		\mathbb{E}\left[T_K^2 \big | X \right] &=  \frac{1}{n^2(n-1)^2} \sum_{i\neq j\neq u\neq v=1}^{n} K(X_i,X_j)K(X_u,X_v) f(X_i)f(X_j)f(X_u)f(X_v)\\
		{} &+ \frac{4}{n^2(n-1)^2} \sum_{i\neq j\neq u=1}^{n} K(X_i,X_j)K(X_i,X_u) \left[f^2(X_i)+\sigma^2 \right]f(X_j)f(X_u)\\
		{} &+ \frac{2}{n^2(n-1)^2} \sum_{i\neq j =1}^{n} K^2(X_i,X_j) \left[f^2(X_i)+\sigma^2 \right] \left[f^2(X_j)+\sigma^2 \right]
		\end{align*}
		Since $\mathbb{E}\left[T_K^2 \right] = \mathbb{E}\left[\mathbb{E}\left[T_K^2 \big | X \right] \right]$, and since $\left(X_1,\cdots,X_n\right)$ are i.i.d with density $\nu$ on $E$, we obtain
		\begin{align*}
		\mathbb{E} \left[T_K^2 \right] &= \frac{(n-2)(n-3)}{n(n-1)}\int_{E^4}^{} K(x,y)K(u,v)f(x)f(y)f(u)f(v)d\nu(x)d\nu(y)d\nu(u)d\nu(v)\\
		{} &+ \frac{4(n-2)}{n(n-1)}\int_{E^3}^{}K(x,y)K(x,u)\left[f^2(x)+\sigma^2 \right]f(y)f(u)d\nu(x)d\nu(y)d\nu(u)\\
		{} &+ \frac{2}{n(n-1)}\int_{E^2}^{}K^2(x,y)\left[f^2(x)+\sigma^2 \right]\left[f^2(y)+\sigma^2 \right]d\nu(x)d\nu(y)\\
		{} &=  \frac{(n-2)(n-3)}{n(n-1)}\varepsilon_K^2 + 4A_K + 2B_K.
		\end{align*}
		Thus
		\begin{equation}\label{emu}
			\mathbb{E}\left[T^2_K\right] - \frac{(n-2)(n-3)}{n(n-1)}\varepsilon_K^2 = 4A_K + 2B_K.
		\end{equation}
		In fact
		\begin{equation*}
		0 <\frac{(n-2)(n-3)}{n(n-1)} <1,\ \forall n>3.
		\end{equation*}
		Then 
		\begin{equation*}
		\mathbb{E}\left[T^2_K\right] - \frac{(n-2)(n-3)}{n(n-1)}\varepsilon_K^2 \geq \mathbb{E}\left[T^2_K\right] - \varepsilon_K^2.
		\end{equation*}
		Replacing (\ref{emu}) into (\ref{P.H1}) we obtain
		\begin{equation}
			\mathbb{P}_{f}\left(| \varepsilon_K-T_K| \geq \lambda\right) \leq \frac{4A_K + 2B_K}{\lambda^2}.
		\end{equation}
		Choosing $\lambda = \sqrt{\frac{16A_K + 8B_K}{\beta}}$, the above inequality leads to
		\begin{equation*} \label{P.ab}
			\mathbb{P}_{f}\left(|\varepsilon_K-T_K| \geq \sqrt{\frac{16A_K + 8B_K}{\beta}}\right) \leq \frac{\beta}{4}.
		\end{equation*}
		This implies
		\begin{equation}\label{P.H1.1}
			\mathbb{P}_{f}\left(T_K \leq \varepsilon_K - \sqrt{\frac{16A_K + 8B_K}{\beta}}\right) \leq \frac{\beta}{4}.
		\end{equation}
		Now we consider the term $\widehat{\sigma}_n^2=\frac{1}{n}\sum_{i=1}^{n}\left(Y_{2i-1}-Y_{2i}\right)^2$. Following to the Cochran's theorem, we consider the orthogonal subspace $W$ of dimension $n/2$. We denote $\left(e_1,\cdots,e_{n/2}\right)$ be an orthogonal basis of $W$, where for all $i=1,\cdots,n/2$, $e^T_i$ is a vetor includes $n$ elements within two values $\{0,1\}$ and its values equal to $1$ at two positions $2i$ and $2i-1$. On the other hand, for $Y=\left(Y_1,\cdots,Y_n\right)$, with $Y_i=f(X_i)+\sigma \epsilon_i$ we have
		\begin{equation*}
		Y=\mathcal{N}_n \left( \left( \begin{array}{c}
		f(X_1)  \\
		\cdots \\
		f(X_n)  \end{array} \right),\sigma^2 I_n \right). 
		\end{equation*}
		Using the Cochran's theorem, we have
		\begin{equation*}
		\norm{\Pi_{W^\bot} Y}^2=\norm{Y-\Pi_W Y}^2 \sim \sigma^2 \chi^2\left(\frac{n}{2},\ \frac{na^2}{2\sigma^2} \right),
		\end{equation*}
		where $a^2 := \frac{1}{n}\sum_{i=1}^{n/2}\left[f\left(\frac{2i-1}{n} \right) -f\left(\frac{2i}{n} \right) \right]^2$ and $\chi^2(k,\lambda)$ denotes a non central Chi-square variable with $k$ degrees of freedom and non centrality parameter $\lambda$.\\
		Moreover,
		\begin{align*}
		\norm{Y-\Pi_W Y}^2 &=\sum_{i=1}^{n/2}\left(Y_{2i-1} - \frac{Y_{2i-1}+Y_2i}{2}\right)^2+\left(Y_{2i} - \frac{Y_{2i-1}+Y_2i}{2}\right)^2\\
		{} &= 2\sum_{i=1}^{n/2} \left(\frac{Y_{2i-1}-Y_{2i}}{2}\right)^2 = \frac{1}{2}\sum_{i=1}^{n/2}\left(Y_{2i-1}-Y_{2i} \right)^2\\
		{} &= \frac{n}{2}\hat{\sigma}_n^2.
		\end{align*}
		Hence
		\begin{equation*}
		\widehat{\sigma}_n^2 \sim \frac{2\sigma^2}{n} \chi^2\left(\frac{n}{2},\ \frac{na^2}{2\sigma^2}  \right).
		\end{equation*}
		Now, we consider the variable $Z\sim \chi^2\left(\frac{n}{2}, \frac{na^2}{2\sigma^2} \right)$. Using Lemma 8.1 in \cite{birge2001alternative}, we have
		\begin{equation*}
		\forall \rho >0,\ \mathbb{P}_{f}\left[Z \geq  \frac{n}{2} + \frac{na^2}{2\sigma^2} + 2\sqrt{\left(\frac{n}{2}+ \frac{na^2}{\sigma^2} \right)\rho } +2\rho  \right] \leq e^{-\rho}.
		\end{equation*}
		This implies
		\begin{equation}\label{P.H1.2.proof}
		\forall \rho >0,\ \mathbb{P}_{f}\left[\hat{\sigma}_n^2  \geq \sigma^2 + a^2 +\frac{4\sigma^2}{n}\sqrt{\left(\frac{n}{2}+\frac{na^2}{\sigma^2} \right)\rho } +\frac{4\rho\sigma^2}{n} \right] \leq e^{-\rho}.
		\end{equation}
		Choosing $\rho=\ln \left(4/\beta\right)$, (\ref{P.H1.2.proof}) leads to
		\begin{equation*}
		\mathbb{P}_{f} \left(\hat{\sigma}_n^2 \geq D_{n,\beta} \right) \leq \frac{\beta}{4},
		\end{equation*}
		where
		\begin{equation*}
		D_{n,\beta}:= \sigma^2 + a^2 +\frac{4\sigma^2}{n}\sqrt{\left(\frac{n}{2} + \frac{na^2}{\sigma^2} \right) \ln \left(\frac{4}{\beta}\right)} +\frac{4\sigma^2}{n}\ln \left(\frac{4}{\beta}\right).
		\end{equation*}
		Thus
		\begin{equation}\label{P.H1.2}
		\mathbb{P}_{f} \left(\frac{1}{\hat{\sigma}_n^2} \leq \frac{1}{D_{n,\beta}} \right) \leq \frac{\beta}{4}.
		\end{equation}
		From (\ref{P.H1.1}) and (\ref{P.H1.2}), we obtain
		\begin{equation}\label{6}
		\mathbb{P}_{f} \left(V_K \leq UV\right)\leq \beta/2,
		\end{equation}
		with
		\begin{align*}
		U&= \varepsilon_K - \sqrt{\frac{16A_K + 8B_K}{\beta}},\\
		V&= \frac{1}{D_{n,\beta}}.
		\end{align*}
		If $q^\alpha_{K,1-\beta/2} \leq UV$ we have $\mathbb{P}_{f}\left( V_K \leq q^\alpha_{K,1-\beta/2}\right) \leq \beta/2$. \\
		Therefore, if
		\begin{equation}
		\varepsilon_K  \geq D_{n,\beta}\  q^\alpha_{K,1-\beta/2}  + \sqrt{\frac{16A_K + 8B_K}{\beta}},
		\end{equation}
		then
		\begin{equation*}
		\mathbb{P}_{f}(\Phi_{K,\alpha} = 0) \leq \beta.
		\end{equation*}
		Let us now give an upper bound for $q^\alpha_{K,1-\beta/2}$. Reasoning conditionally on $X$, we recognize in $\frac{1}{n(n-1)}\sum_{i\neq j}^{n}K\left(X_i,X_j\right)\epsilon_i \epsilon_j:= T_K^{\left(0\right)}$ be a Gaussian chaos, as defined by \citeauthor{de2012decoupling} \citeyearpar{de2012decoupling}, of the form $Z= \sum_{i \neq i^{'}}^{} x_{i,i^{'}} \epsilon_i \epsilon_{i^{'}}$, where $x_{i,i^{'}}$'s are some real deterministic numbers and $(\epsilon_i)_i$ is a sequence of i.i.d Gaussian variables.  Corollary 3.26 of \cite{de2012decoupling} states that there exists some absolute constant $\kappa>0$ such that if $\gamma^2 = \mathbb{E}[Z^2] =  \sum_{i \neq i^{'}}^{} x_{i,i^{'}}^2$.
		Then
		\begin{equation*}
		\mathbb{E}\left[\exp \left(\frac{|Z|}{\kappa \gamma}\right)\right] \leq 2.
		\end{equation*}
		Hence by Markov's inequality,
		\begin{equation}\label{TK.H0}
		\mathbb{P}\left(|Z| \geq \kappa\gamma \ln \left(\frac{2}{\alpha}\right)\right)\leq \alpha.
		\end{equation}
		Applying the result (\ref{TK.H0}) for $T_K^{\left(0\right)}$ with 
		\begin{equation}
		\gamma^2 = \frac{\sigma^4}{n^2(n-1)^2}\sum_{i\neq j =1}^{n} K^2_{ij},
		\end{equation}
		we have
		\begin{equation}\label{P.H0.1}
		\mathbb{P}_{(H_0)}\left(T_K \geq \frac{\kappa \sigma^2 }{n(n-1)} \ln \left(\frac{4}{\alpha}\right) \sqrt{\sum_{i\neq j =1}^{n} K^2_{ij}} \right) \leq \frac{\alpha}{2}.
		\end{equation}
		On the other hand, we have, under $H_0$, $\hat{\sigma}_n^2=\hat{\sigma}^2_{n,\epsilon}$ where
		\begin{equation*}
		\widehat{\sigma}_{n,\epsilon}^2 = \frac{\sigma^2}{n} \sum_{i=1}^{n/2} \left(\epsilon^{'}_{2i-1} - \epsilon^{'}_{2i}\right)^2= \frac{2\sigma^2}{n}\sum_{i=1}^{n/2}\left(\frac{\epsilon^{'}_{2i-1} - \epsilon^{'}_{2i}}{\sqrt{2}}\right)^2 \sim \frac{2\sigma^2}{n} \chi_{(n/2)}.
		\end{equation*}
		Since the variables $Z_i=\epsilon^{'}_{2i-1} - \epsilon^{'}_{2i}$, $i=1,\cdots,n/2$ are i.i.d standard Gaussian variables.
		Using the Lemma 8.1 in \cite{birge2001alternative}, we obtain
		\begin{equation}\label{P-e.x}
		\forall x>0,\ \mathbb{P}_{(H_0)} \left(\frac{2\sigma^2}{n} \chi_{(n/2)} \leq \sigma^2 -\frac{4\sigma^2}{\sqrt{2n}}\sqrt{x} \right) \leq e^{-x}.
		\end{equation}
		Choosing $x=\ln(2/\alpha)$, (\ref{P-e.x}) leads to
		\begin{equation}\label{P.H0.2}
		\mathbb{P}_{(H_0)}\left(\hat{\sigma}_n^2 \leq \sigma^2 -\frac{4\sigma^2}{\sqrt{2n}}\sqrt{\ln \left(\frac{2}{\alpha}\right)} \right) \leq \frac{\alpha}{2}.
		\end{equation}
		Moreover, we have
		\begin{align*}
		\mathbb{P}_{(H_0)} \left(\frac{T_K}{\hat{\sigma}_n^2} \geq \frac{\frac{\kappa \sigma^2 }{n(n-1)} \ln \left(\frac{2}{\alpha}\right) \sqrt{\sum_{i\neq j =1}^{n} K^2_{ij}}}{\sigma^2\left(1-\frac{2\sqrt{2}}{\sqrt{n}}\sqrt{\ln \left(\frac{2}{\alpha}\right)}\right) } \right) &\leq \mathbb{P}_{(H_0)}\left(T_K \geq \frac{\kappa \sigma^2 }{n(n-1)} \ln \left(\frac{2}{\alpha}\right) \sqrt{\sum_{i\neq j =1}^{n} K^2_{ij}} \right) \\
		{} &+ \mathbb{P}_{(H_0)} \left(\frac{1}{\hat{\sigma}_n^2} \geq \frac{1}{\sigma^2\left(1-\frac{2\sqrt{2}}{\sqrt{n}}\sqrt{\ln \left(\frac{2}{\alpha}\right)}\right)} \right).
		\end{align*}
		From (\ref{P.H0.1}) and (\ref{P.H0.2}), we obtain
		\begin{equation}
		\mathbb{P}_{(H_0)} \left(\frac{T_K}{\hat{\sigma}_n^2} \geq \frac{\frac{\kappa \sigma^2 }{n(n-1)} \ln \left(\frac{2}{\alpha}\right) \sqrt{\sum_{i\neq j =1}^{n} K^2_{ij}}}{\sigma^2\left(1-\frac{2\sqrt{2}}{\sqrt{n}}\sqrt{\ln \left(\frac{2}{\alpha}\right)}\right) } \right) \leq \frac{\alpha}{2} +\frac{\alpha}{2} =\alpha.
		\end{equation}
		This implies
		\begin{equation*}
		\mathbb{P}_{(H_0)} \left(V_K \geq \frac{\frac{\kappa \sigma^2 }{n(n-1)} \ln \left(\frac{2}{\alpha}\right) \sqrt{\sum_{i\neq j =1}^{n} K^2_{ij}}}{\sigma^2\left(1-\frac{2\sqrt{2}}{\sqrt{n}}\sqrt{\ln \left(\frac{2}{\alpha}\right)}\right) } \right) \leq \frac{\alpha}{2} +\frac{\alpha}{2} =\alpha.
		\end{equation*}
		Thus the $(1-\alpha)$ quantile of $V_K^{\left(0\right)}$ conditionally on $X$ satisfies
		\begin{equation}\label{quant}
		q_{K,1-\alpha}^{(X)} \leq \frac{\frac{\kappa \sigma^2 }{n(n-1)} \ln \left(\frac{2}{\alpha}\right) \sqrt{\sum_{i\neq j =1}^{n} K^2_{ij}}}{\sigma^2\left(1-\frac{2\sqrt{2}}{\sqrt{n}}\sqrt{\ln \left(\frac{2}{\alpha}\right)}\right) }.
		\end{equation}
		Taking $n\geq 32\ln \left(\frac{2}{\alpha}\right)$, so $\sqrt{n} \geq 4\sqrt{2}\sqrt{\ln \left(\frac{2}{\alpha}\right)} $, (\ref{quant}) returns to
		\begin{equation}
		q_{K,1-\alpha}^{(X)} \leq \frac{2\kappa}{n(n-1)} \ln \left(\frac{2}{\alpha}\right) \sqrt{\sum_{i\neq j =1}^{n} K^2_{ij}}.
		\end{equation}
		Hence $q^\alpha_{K,1-\beta/2}$ is upper bounded by the $(1-\beta/2)$ quantile of $\frac{2\kappa}{\sqrt{n(n-1)}} \ln \left(\frac{2}{\alpha}\right) \sqrt{\frac{1}{n(n-1)} \sum_{i\neq j =1}^{n} K^2_{ij}}$.\\
		We define 
		\begin{equation*}
		U_n = \frac{1}{n(n-1)} \sum_{i\neq j =1}^{n} K^2_{ij}.
		\end{equation*}
		We use Markov's inequality again for the nonnegative random variable $U_n$, we obtain for any $\delta>0$
		\begin{equation}\label{u.nP}
		\mathbb{P}_{f} \left( U_n >\delta \right) \leq \frac{\mathbb{E}(U_n) }{\delta}.
		\end{equation}
		We have
		\begin{equation}
		\mathbb{E}\left[U_n\right] = \int_{E^2}^{}K^2(x,y)d\nu(x)d\nu(y).
		\end{equation}
		Choosing $\delta=2\int_{E^2}^{}K^2(x,y)d\nu(x)d\nu(y)/\beta$, (\ref{u.nP}) returns to 
		\begin{equation*}
		\mathbb{P}_{f} \left( U_n >\frac{2\int_{E^2}^{}K^2(x,y)d\nu(x)d\nu(y)}{\beta} \right) \leq \frac{\beta}{2},
		\end{equation*}
		and
		\begin{equation*}
		q^\alpha_{K,1-\beta/2} \leq \frac{2\kappa}{\sqrt{n(n-1)}} \ln \left(\frac{2}{\alpha}\right) \sqrt{\frac{2\int_{E^2}^{}K^2(x,y)d\nu(x)d\nu(y)}{\beta} },
		\end{equation*}
		which concludes the proof.
		\section{Proof of Theorem \ref{thmnormdiff}}
		For all symmetric kernel function $K$, we have
		\begin{equation}\label{eps}
			\mathbb{E}\left(T_K\right)=\langle K[f],f\rangle = \frac{1}{2}\left(\| f\|^2 +\|K[f] \|^2 -\|f-K[f] \|^2 \right).
		\end{equation}
		On the other hand
		\begin{equation*}
		A_K \leq \frac{\left( \|f\|^2_\infty+\sigma^2 \right)\norm{K[f]}^2}{n} .
		\end{equation*}
		Let $C_K$ be an upper bound for $\int_{E^2}^{}K^2(x,y)d\nu(x)d\nu(y)$, we have
		\begin{equation*}
		B_K \leq \frac{\left( \|f\|^2_\infty+\sigma^2 \right)^2C_K}{n(n-1)}
		\end{equation*}
		From Proposition \ref{2nderror}, the bounds for $A_K$ and $B_K$ and the inequality $\sqrt{a+b}\leq \sqrt{a}+\sqrt{b}$ for all $a\geq 0,b\geq 0$, we deduce that $\mathbb{P}_{f}\left(\Phi_{K,\alpha}=0 \right)\leq \beta$ as soon as,
		\begin{align*}
	\| f\|^2 +\|K[f] \|^2 -\|f-K[f] \|^2  &\geq 4\sqrt{2}\kappa D_{n,\beta}  \ln \left(\frac{2}{\alpha}\right) \sqrt{\frac{C_K}{n(n-1)\beta}} + 4\sqrt{2} \left( \norm{f}^2_\infty+\sigma^2\right)  \sqrt{\frac{C_K}{n(n-1)\beta}}\\
		{} &+ 8\norm{K[f]}\sqrt{\frac{\left( \|f\|^2_\infty+\sigma^2 \right)}{n\beta}}.
		\end{align*}
		By using the elementary inequality $2cd \leq c^2+d^2$ with $c=\norm{K[f]}$ and $d=4\sqrt{\frac{\left( \|f\|^2_\infty+\sigma^2 \right)}{n\beta}}$ in the right hang side of the above condition, the above condition holds if
		\begin{align*}
		\| f\|^2  &\geq \|f-K[f] \|^2 + \frac{16\left( \|f\|^2_\infty+\sigma^2 \right)}{n\beta}\\
		{} &+ \frac{4\sqrt{C_K}}{\sqrt{n(n-1)\beta}} \left(\kappa D_{n,\beta}\ln \left(\frac{2}{\alpha}\right) +\sqrt{2}\left( \|f\|^2_\infty+\sigma^2 \right) \right).
		\end{align*}

			\section{Proof of Corollary \ref{col.f.gaussian} and \ref{col.gauss}.}
			Under the hypothesis of corollary \ref{col.f.gaussian},
			\begin{equation*}
			K(x,y)=\sum_{\lambda \in \Lambda}^{} \phi_\lambda(x) \phi_\lambda(y),
			\end{equation*}
			and the linear space $S$ generated by the functions $\left( \phi_\lambda,\lambda \in \Lambda \right)$ is of dimension $D$. Hence, we have
			\begin{equation*}
			\int_{[0,1]^2}^{}\left(\sum_{\lambda \in \Lambda}^{} \phi_\lambda(x) \phi_\lambda(y) \right)^2d\nu(x)d\nu(y)\leq D.
			\end{equation*}
			Thus, we can take $C_K=D$.\\
			Second, under choice of the Gaussian kernel defined by (\ref{13}), we recall that
			\begin{equation}
			K(x,y)=\frac{1}{h}k\left(\frac{x-y}{h}\right),\quad \text{for } \left(x,y\right)\in \mathbb{R}^2
			\end{equation}
			where $k(u)=\frac{1}{\sqrt{2\pi}}\exp \left(-u^2/2 \right),\ \text{for all} \ u\in \mathbb{R}$ and $h$ is a positive bandwidth.\\
			We have
			\begin{align*}
			\int_{\mathbb{R}^2}^{}K^2(x,y)\nu(x)\nu(y)d(x)d(y) &= \frac{1}{2\pi h^2}\int_{\mathbb{R}^2}^{}e^{-\frac{(x-y)^2}{2h^2}}\nu(x)\nu(y)d(x)d(y)\\
			&\leq \frac{1}{2\pi h^2}\int_{\mathbb{R}^2}^{}e^{-\frac{u^2}{2}}h du \nu(x)\nu(x-uh)dx\\
			&\leq \frac{\norm{\nu}_\infty}{h\sqrt{2\pi}} \left(\frac{1}{\sqrt{2\pi}}\int_{\mathbb{R}}^{}e^{-\frac{u^2}{2}} du \right)\left(\int_{\mathbb{R}}^{}\nu(x)dx \right)\\
			&\leq \frac{\norm{\nu}_\infty}{h\sqrt{2\pi}}.
			\end{align*}
			Hence, we can choose $C_K=\frac{\norm{\nu}_\infty}{h\sqrt{2\pi}}$.
			
		\section{Proof of Proposition \ref{prop.project}}
				For all $J\geq 0$, we set $D=2^J$ be the dimension of $S_J$. Let us assume that $f \in \mathcal{B}_{2,\infty}^\delta (R)$, it implies 
				\begin{equation}
				\|f-\Pi_{S_J}(f) \|^2 \leq R^2 2^{-2J\delta}.
				\end{equation}
				We obtain from Corollary \ref{col.f.gaussian} that there exists
				\begin{equation*}
				C(\alpha,\beta,\sigma,\|f\|_\infty) >0,
				\end{equation*}
				such that $\mathbb{P}_{f} \left(\Phi_{K,\alpha}=0 \right)\leq \beta$ if
				\begin{equation}\label{normproj1}
				\|f\|^2  \geq C(\alpha,\beta,\sigma,R,\|f\|_\infty)  \left[ 2^{-2J\delta}+ \frac{2^{J/2}}{n} \right].
				\end{equation}
				In this case, we see that the right hand side of (\ref{normproj}) reproduces a bias-variance decomposition close to the bias-variance decomposition for projection estimators, with the bias term $R^2 2^{-2J\delta}$ and the variance term $2^{J/2}/n$. The optimal choice of $J$ satisfies
				\begin{equation*}
				R^2 2^{-2J\delta} = \frac{2^{J/2}}{n}.
				\end{equation*}
				Thus, we obtain the optimal choice $J^{*}$,
				\begin{equation*}
				J^{*} = \left[\log_2\left(n^{2/(1+4\delta)} \right)\right],
				\end{equation*}	
				leading to the desired result.
	
	\section{Proof of Proposition \ref{prop.gauss}}
		Considering (\ref{normgaus}), we mainly have to find a sharp upper bound for $\|f-k_h*f \|^2$ when $f\in \mathcal{S}_\delta(R)$. Plancherel's theorem gives that when $f \in \mathbb{L}^1(\mathbb{R})\cap \mathbb{L}^2(\mathbb{R})$,
		\begin{align*}
		(2\pi)\|f-k_h*f \|^2 &= \norm{\left(1-\widehat{k}_h \right) \left(\widehat{f} \right) }^2 \\
		{} &= \int_{\mathbb{R}}^{}\bigg |1- \widehat{k}(2^{-l}u) \bigg|^2(u) {\hat{f}}^2(u) du.
		\end{align*}
		We assume that $\norm{\hat{k}}_\infty <+\infty$ and
		\begin{equation*}
		\text{Ess} \sup_{u\in \mathbb{R}\setminus \{0\} } \frac{|1 - \widehat{k}(u)| }{|u|} \leq C,
		\end{equation*}
		for some $C >0$. There also exists some constant $C(\delta)>0$ such that
		\begin{equation*}
		\text{Ess} \sup_{u\in \mathbb{R}\setminus \{0\} } \frac{|1 - \widehat{k}(u)| }{|u|^\delta} \leq C(\delta).
		\end{equation*}
		Then
		\begin{equation*}
		\|f-k_h*f \|^2 \leq \frac{C(\delta)}{2\pi} \int_{\mathbb{R}}^{}|2^{-l}u|^{2\delta} \hat{f}^2(u)du
		\end{equation*}
		and since $f \in \mathcal{S}_\delta(R)$,
		\begin{equation*}
		\|f-k_h*f \|^2 \leq 2^{-2\delta l}C(\delta) R^2.
		\end{equation*}
		We obtain from corollary \ref{col.gauss} that there exists 
		\begin{equation*}
		C(\alpha,\beta,\sigma,R,\|f\|_\infty) >0,
		\end{equation*}
		such that $\mathbb{P}_{f} \left(\Phi_{K,\alpha}=0 \right)\leq \beta$ if
		\begin{equation}\label{normgauss1}
		\|f\|^2  \geq C(\alpha,\beta,\sigma,R,\|f\|_\infty)  \left[2^{-2\delta l} + \frac{2^{l/2}}{n} \right].
		\end{equation}
		In this case, we see that the right hand side of (\ref{normgauss1}) reproduces a bias-variance decomposition with the bias term $2^{-2\delta l}$ and the variance term $2^{l/2}/n$. The optimal choice of $l$ satisfies
		\begin{equation*}
		2^{-2\delta l} = \frac{2^{l/2}}{n}.
		\end{equation*}
		Thus, we obtain the optimal choice $l^{*}$ as follows.
		\begin{equation*}
		l^{*} = \left[\log_2\left(n^{2/(1+4\delta)} \right) \right],
		\end{equation*}
		leading to the desired result.
		
	\section{Proof of Theorem \ref{them4.1them}, Corollary \ref{colaggreProj} and \ref{gauss.aggre}.}		
	We have 
	\begin{align*}
	\mathbb{P}_{(H_0)} \left(\Phi_{\alpha}=1 \right) &=\mathbb{P}\left[\sup_{m\in \mathcal{M}}\left(V_{K_m}-q^{(X)}_{m,1-u_\alpha^{(X)} e^{-w_m}} \right)>0 \right]\\
	{} &= \mathbb{E}\left[\mathbb{E} \left[ \mathbbm{1}\left\lbrace \sup_{m\in \mathcal{M}}\left(V_{K_m}-q^{(X)}_{m,1-u_\alpha^{(X)} e^{-w_m}} \right)>0 \right\rbrace \bigg | X \right] \right].
	\end{align*}
	We have,
	\begin{equation*}
	\mathbb{E} \left[ \mathbbm{1}\left\lbrace \sup_{m\in \mathcal{M}}\left(V_{K_m}-q^{(X)}_{m,1-u_\alpha^{(X)} e^{-w_m}} \right)>0 \right\rbrace \bigg | X \right] \leq \alpha,
	\end{equation*}
	by definition of $u_\alpha^{(X)} $, which implies that $\mathbb{P}_{(H_0)} \left(\Phi_{\alpha}=1 \right) \leq \alpha$.\\
	On the other hand, we know that $u_\alpha^{(X)} \geq \alpha$. Setting $\alpha_m = \alpha e^{-w_m}$, we have
	\begin{align*}
	\mathbb{P} &{} \left(\exists m\in \mathcal{M}, V_{K_m} > q^{(X)}_{K_m,1-u_\alpha^{(X)} e^{-w_m}} \right) \geq \mathbb{P}\left(\exists m\in \mathcal{M}, V_{K_m} > q^{(X)}_{K_m,1-\alpha_m} \right)\\
	{} &\geq 1- \mathbb{P}\left(\forall m \in \mathcal{M}, V_{K_m} \leq q^{(X)}_{K_m,1-\alpha_m} \right) \geq 1-\inf_{m\in \mathcal{M}} \mathbb{P}\left(V_{K_m}\leq q^{(X)}_{K_m,1-\alpha_m} \right)\\
	{} &\geq 1-\beta,
	\end{align*}
	as soon as there exists $m$ in $\mathcal{M}$ such that
	\begin{equation*}
	\mathbb{P}\left(V_{K_m} \leq q^{(X)}_{K_m,1-\alpha_m}\right) \leq \beta.
	\end{equation*}
	We can now apply Corollary \ref{col.f.gaussian} and \ref{col.gauss} with $\alpha_m = \alpha e^{-w_m}$, so we replace $\ln(2/\alpha)$ by $\left(\ln(2/\alpha)+w_m \right)$ for desired results in Corollary \ref{colaggreProj} and \ref{gauss.aggre}.
	\section{Proof of Corollary \ref{corol4.3}.}
		Considering (\ref{aggrePrker}), we aim to find a sharp upper bound for the right hand side of the inequality when $f\in \mathcal{B}_{2,\infty}^\delta (R,R^{'})$. Let us assume that $f\in \mathcal{B}_{2,\infty}^\delta (R,R^{'})$. Then $f\in\mathcal{B}_{2,\infty}^\delta (R)$, we have
		\begin{equation*}
		\|f-\Pi_{S_J}(f) \|^2 \leq C(\delta)R^22^{-2\delta J},
		\end{equation*}
		and
		\begin{equation*}
		C(\alpha,\beta,\sigma,R,\|f\|_\infty) \leq C(\alpha,\beta,\sigma,R,R^{'}).
		\end{equation*}
		Hence (\ref{aggrePrker}) can be upper bounded by
		\begin{equation*}
		C(\alpha,\beta,\sigma,R,R^{'}) \inf_{J\in \mathcal{M}_{\bar{J}}} \left\lbrace 2^{-2\delta J}+\ln (2+J)\frac{2^{J/2}}{n} \right\rbrace .
		\end{equation*}
		Taking
		\begin{equation*}
		J^{**}=\left[ \log_2\left(\left(\frac{n}{\ln \ln(n)} \right)^{\frac{2}{4\delta+1}} \right) \right],
		\end{equation*}
		\begin{align*}
		C(\alpha,\beta,\sigma,R,R^{'}) &\quad  \inf_{J\in \mathcal{M}_{\bar{J}}} \left\lbrace 2^{-2\delta J}+\ln (2+J)\frac{2^{J/2}}{n} \right\rbrace \\
		{} &\leq C(\alpha,\beta,\sigma,R,R^{'})  \left\lbrace 2^{-2\delta J^{**}}+\ln (2+J^{**})\frac{2^{J^{**}/2}}{n}   \right\rbrace \\
		{} &\leq C(\alpha,\beta,\sigma,R,R^{'}) \left(\frac{n}{\ln \ln(n)} \right)^{-\frac{4\delta}{4\delta+1}}.
		\end{align*}
		That leads to $\mathbb{P}_{f}\left(\Phi_\alpha^{(1)}=0 \right) \leq \beta$ if
		\begin{equation*}
		\norm{f}^2 \geq C(\alpha,\beta,\sigma,R,R^{'})\left(\frac{\ln \ln(n)}{n}\right)^{\frac{4\delta}{4\delta+1}}.
		\end{equation*}
	
		\section{Proof of Corollary \ref{gauss.aggre}.}

		Considering (\ref{gaussiAggre}), we aim to find a sharp upper bound for the right hand side of the inequality when $f\in  \mathcal{S}_\delta(R,R^{'})$. Let us assume that $f\in  \mathcal{S}_\delta(R,R^{'})$. Similarly, with regards to the proof of Proposition \ref{prop.gauss}, we have 
		\begin{equation*}
		\|f-k_l*f \|^2 \leq 2^{-2\delta l}C(\delta)R^2,
		\end{equation*}
		and 
		\begin{equation*}
		C(\alpha,\beta,\sigma,\norm{f}_\infty) \leq C(\alpha,\beta,\sigma,R^{'}).
		\end{equation*}
		Hence (\ref{gaussiAggre}) can be upper bounded by
		\begin{equation*}
		C(\alpha,\beta,\sigma,R,R^{'}) \inf_{l\in \mathcal{L}} \left\lbrace  2^{-2\delta l} + \frac{w_l}{n\sqrt{2^{-l}}} \right\rbrace .
		\end{equation*}
		Choosing 
		\begin{equation*}
		l^{**}=\left[ \log_2\left(\left(\frac{n}{\ln \ln(n)}\right)^{\frac{2}{1+4\delta}} \right) \right],
		\end{equation*}
		\begin{align*}
		C(\alpha,\beta,\sigma,R,R^{'}) &\ \inf_{l\in \mathcal{L}} \left\lbrace  2^{-2\delta l} + \frac{w_l}{n\sqrt{2^{-l}}} \right\rbrace \\
		{} &\leq C(\alpha,\beta,\sigma,R,R^{'}) \left\lbrace  2^{-2\delta l^{**}} + \frac{w_{l^{**}}}{n\sqrt{2^{-l^{**}}}} \right\rbrace \\
		{} &\leq C(\alpha,\beta,\sigma,R,R^{'}) \left( \frac{n}{\ln\ln(n)}\right)^{-\frac{4\delta}{4\delta+1}}.
		\end{align*}
		That leads to $\mathbb{P}_{f}\left(\Phi_\alpha^{(1)}=0 \right) \leq \beta$ if
		\begin{equation*}
		\norm{f}^2 \geq C(\alpha,\beta,\sigma,R,R^{'})\left(\frac{\ln \ln(n)}{n}\right)^{\frac{4\delta}{4\delta+1}}.
		\end{equation*}

	\end{appendices}
	\section*{Acknowledgement}
	I gratefully thank to Professor B\'eatrice Laurent of Institut National des Sciences Appliqu\'ees de Toulouse and Professor Jean-Michel Loubes of Institut de Math\'ematiques de Toulouse for supporting me in the best ideas and comments. 

	\nocite{*}
	\bibliographystyle{apalike} 
	\bibliography{twosample.bib}

\begin{thebibliography}{}

\bibitem[Bachoc et~al., 2017]{bachoc2017gaussian}
Bachoc, F., Gamboa, F., Loubes, J.-M., and Venet, N. (2017).
\newblock A gaussian process regression model for distribution inputs.
\newblock {\em IEEE Transactions on Information Theory}.

\bibitem[Baraud et~al., 2002]{baraud2002non}
Baraud, Y. et~al. (2002).
\newblock Non-asymptotic minimax rates of testing in signal detection.
\newblock {\em Bernoulli}, 8(5):577--606.

\bibitem[Baraud et~al., 2003]{baraud2003adaptive}
Baraud, Y., Huet, S., Laurent, B., et~al. (2003).
\newblock Adaptive tests of linear hypotheses by model selection.
\newblock {\em The Annals of Statistics}, 31(1):225--251.

\bibitem[Birg{\'e}, 2001]{birge2001alternative}
Birg{\'e}, L. (2001).
\newblock An alternative point of view on lepski's method.
\newblock {\em Lecture Notes-Monograph Series}, pages 113--133.

\bibitem[Butucea and Tribouley, 2006]{butucea2006nonparametric}
Butucea, C. and Tribouley, K. (2006).
\newblock Nonparametric homogeneity tests.
\newblock {\em Journal of statistical planning and inference}, 136(3):597--639.

\bibitem[Castillo et~al., 2006]{castillo2006exact}
Castillo, I., L{\'e}vy-Leduc, C., and Matias, C. (2006).
\newblock Exact adaptive estimation of the shape of a periodic function with
  unknown period corrupted by white noise.
\newblock {\em Mathematical methods of statistics}, 15(2):146--175.

\bibitem[De~la Pena and Gin{\'e}, 2012]{de2012decoupling}
De~la Pena, V. and Gin{\'e}, E. (2012).
\newblock {\em Decoupling: from dependence to independence}.
\newblock Springer Science \& Business Media.

\bibitem[Delgado, 1992]{delgado1992testing}
Delgado, M.~A. (1992).
\newblock Testing the equality of nonparametric regression curves.

\bibitem[Eubank and LaRiccia, 1993]{eubank1993testing}
Eubank, R. and LaRiccia, V. (1993).
\newblock Testing for no effect in nonparametric regression.
\newblock {\em Journal of statistical planning and inference}, 36(1):1--14.

\bibitem[Fromont et~al., 2011]{fromont2011adaptive}
Fromont, M., Laurent, B., Reynaud-Bouret, P., et~al. (2011).
\newblock Adaptive tests of homogeneity for a poisson process.
\newblock In {\em Annales de l'Institut Henri Poincar{\'e}, Probabilit{\'e}s et
  Statistiques}, volume~47, pages 176--213. Institut Henri Poincar{\'e}.

\bibitem[Fromont et~al., 2013]{fromont2013two}
Fromont, M., Laurent, B., Reynaud-Bouret, P., et~al. (2013).
\newblock The two-sample problem for poisson processes: Adaptive tests with a
  nonasymptotic wild bootstrap approach.
\newblock {\em The Annals of Statistics}, 41(3):1431--1461.

\bibitem[Fromont et~al., 2012]{fromont2012kernels}
Fromont, M., Lerasle, M., Reynaud-Bouret, P., et~al. (2012).
\newblock Kernels based tests with non-asymptotic bootstrap approaches for
  two-sample problems.
\newblock In {\em Conference on Learning Theory}, pages 23--1.

\bibitem[Fromont and L{\'e}vy-Leduc, 2006]{fromont2006adaptive}
Fromont, M. and L{\'e}vy-Leduc, C. (2006).
\newblock Adaptive tests for periodic signal detection with applications to
  laser vibrometry.
\newblock {\em ESAIM: Probability and Statistics}, 10:46--75.

\bibitem[Gretton et~al., 2007]{gretton2007kernel}
Gretton, A., Borgwardt, K.~M., Rasch, M., Sch{\"o}lkopf, B., and Smola, A.~J.
  (2007).
\newblock A kernel method for the two-sample-problem.
\newblock In {\em Advances in neural information processing systems}, pages
  513--520.

\bibitem[Hall and Hart, 1990]{hall1990bootstrap}
Hall, P. and Hart, J.~D. (1990).
\newblock Bootstrap test for difference between means in nonparametric
  regression.
\newblock {\em Journal of the American Statistical Association},
  85(412):1039--1049.

\bibitem[Hardle and Marron, 1990]{hardle1990semiparametric}
Hardle, W. and Marron, J.~S. (1990).
\newblock Semiparametric comparison of regression curves.
\newblock {\em The Annals of Statistics}, pages 63--89.

\bibitem[Huskova and Janssen, 1993]{huskova1993consistency}
Huskova, M. and Janssen, P. (1993).
\newblock Consistency of the generalized bootstrap for degenerate u-statistics.
\newblock {\em The Annals of Statistics}, pages 1811--1823.

\bibitem[Ingster, 1993]{ingster1993asymptotically}
Ingster, Y.~I. (1993).
\newblock Asymptotically minimax hypothesis testing for nonparametric
  alternatives. i, ii, iii.
\newblock {\em Math. Methods Statist}, 2(2):85--114.

\bibitem[King et~al., 1991]{king1991testing}
King, E., Hart, J.~D., and Wehrly, T.~E. (1991).
\newblock Testing the equality of two regression curves using linear smoothers.
\newblock {\em Statistics \& Probability Letters}, 12(3):239--247.

\bibitem[King, 1988]{king1988test}
King, E.~C. (1988).
\newblock {\em A test for the equality of two regression curves based on kernel
  smoothers}.
\newblock PhD thesis, Texas A \& M University.

\bibitem[Peyr{\'e} et~al., 2016]{peyre2016gromov}
Peyr{\'e}, G., Cuturi, M., and Solomon, J. (2016).
\newblock Gromov-{W}asserstein averaging of kernel and distance matrices.
\newblock In {\em ICML 2016}.

\bibitem[Romano and Wolf, 2005]{romano2005exact}
Romano, J.~P. and Wolf, M. (2005).
\newblock Exact and approximate stepdown methods for multiple hypothesis
  testing.
\newblock {\em Journal of the American Statistical Association},
  100(469):94--108.

\bibitem[Spokoiny et~al., 1996]{spokoiny1996adaptive}
Spokoiny, V.~G. et~al. (1996).
\newblock Adaptive hypothesis testing using wavelets.
\newblock {\em The Annals of Statistics}, 24(6):2477--2498.

\bibitem[Tsybakov, 2008]{Tsybakov:2008:INE:1522486}
Tsybakov, A.~B. (2008).
\newblock {\em Introduction to Nonparametric Estimation}.
\newblock Springer Publishing Company, Incorporated, 1st edition.

\bibitem[Villani, 2008]{villani2008optimal}
Villani, C. (2008).
\newblock {\em Optimal transport: old and new}, volume 338.
\newblock Springer Science \& Business Media.

\bibitem[Whitt, 1976]{whitt1976bivariate}
Whitt, W. (1976).
\newblock Bivariate distributions with given marginals.
\newblock {\em The Annals of statistics}, pages 1280--1289.

\end{thebibliography}

\end{document}